\documentclass[preprint]{imsart}

\RequirePackage[OT1]{fontenc}

\usepackage{amsthm,amsmath}

\usepackage{amssymb}
\usepackage{amsfonts}

\usepackage{mathrsfs}

\usepackage{epsfig}
\usepackage{graphicx}
\usepackage{subfig}

\usepackage{color}

\usepackage{makecell}
\usepackage{multirow}

\usepackage{url}

\usepackage[toc,page]{appendix}

\RequirePackage[numbers]{natbib}
\RequirePackage[colorlinks,citecolor=blue,urlcolor=blue]{hyperref}

\doi{10.1214/154957804100000000}
\pubyear{0000}
\volume{0}
\firstpage{0}
\lastpage{0}
\arxiv{}

\startlocaldefs

\newcommand{\ton}{\underset{n \to \infty}{\longrightarrow}}

\newcommand{\xb}{\mathbf{x}}
\newcommand{\cb}{\mathbf{c}}
\newcommand{\XX}{\mathscr{X}}

\newcommand{\rr}{\mathscr{L}}

\newcommand{\NN}{\mathbb{N}}

\newcommand{\ee}{\mathbb{E}}
\newcommand{\al}{\boldsymbol\alpha}
\newcommand{\be}{\boldsymbol\beta}
\newcommand{\ub}{\mathcal{U}} 
\newcommand{\vb}{\mathcal{V}} 

\newcommand{\R}{\mathbb{R}}

\newcommand{\E}{\mathbb{E}}
\newcommand{\V}{\mathbb{V}}

\newcommand{\bPsi}{\boldsymbol{\Psi}} 
\newcommand{\etab}{\boldsymbol{\eta}} 

\newcommand{\cov}{\mathrm{Cov}}

\newcommand{\UU}{\mathcal{A}} 
\newcommand{\unUU}{\mathcal{B}} 

\newcommand{\unsim}{\mathord{\sim}}

\newcommand{\pnu}{V}

\newcommand{\T}{\mathbb{L}}
\newcommand{\TT}{\mathbb{T}}

\newcommand{\la}{\langle}
\newcommand{\ra}{\rangle}

\newcommand{\eq}{\triangleq}
\newcommand{\brho}{\boldsymbol\rho}

\newcommand{\bchi}{\Omega^{+}}

\newcommand{\bz}{\mathbf{z}}

\newcommand{\des}{\mathcal{D}}

\DeclareMathOperator*{\argmin}{arg\,min}

\DeclareMathOperator{\sign}{sign}

\newcommand{\mer}{\mathscr{E}} 
\newcommand{\ris}{\mathcal{R}} 

\newcommand{\eps}{\varepsilon}

\newcommand*{\QEDA}{\hfill\ensuremath{\blacksquare}}%

\let\emptyset\varnothing
\newtheorem{definition}{Definition}
\newtheorem{remark}{Remark}
\newtheorem{theorem}{Theorem}
\newtheorem{corollary}{Corollary}
\newtheorem{example}{Example}
\newtheorem{Condition}{Condition}
\newtheorem{lemma}{Lemma}

\theoremstyle{plain}

\numberwithin{equation}{section}
\graphicspath{{pics/}}

\endlocaldefs

\begin{document}
\sloppy 

\begin{frontmatter}

\title{Risk of estimators for Sobol' sensitivity indices based on metamodels}

\runtitle{Risk of estimators for Sobol' indices}


\author{\fnms{Ivan} \snm{Panin}\thanksref{t1}\ead[label=e1]{panin@phystech.edu}}
\address{Skolkovo Institute of Science and Technology\\ 
Bolshoy Boulevard 30, bld. 1, 
Moscow, Russia 121205 \\
Kharkevich Institute for Information Transmission Problems\\
Bolshoy Karetny per. 19, bld. 1, Moscow, Russia 127051\\ \printead{e1}}

\thankstext{t1}{The author would like to thank Maxim Panov for his helpful comments and support.}

\runauthor{Ivan Panin}

\begin{abstract}
Sobol' sensitivity indices allow to quantify the respective effects of random input variables and their combinations on the variance of mathematical model output. We focus on the~problem of Sobol' indices estimation via a metamodeling approach where we replace the true mathematical model with a sample-based approximation to compute sensitivity indices. We propose a new method for indices quality control and obtain asymptotic and non-asymptotic risk bounds for Sobol' indices estimates based on a general class of metamodels. Our analysis is closely connected with the~problem of nonparametric function fitting using the orthogonal system of functions in the random design setting. It considers the relation between the metamodel quality and the error of the corresponding estimator for Sobol' indices and shows the possibility of fast convergence rates in the case of noiseless observations. The theoretical results are complemented with numerical experiments for the approximations based on multivariate Legendre and Trigonometric polynomials.
\end{abstract}



\begin{keyword}[class=MSC]
\kwd[Primary ]{62J10}
\kwd{62J05}
\kwd[; secondary ]{65T40}
\end{keyword}


\begin{keyword}
\kwd{Global Sensitivity Analysis}
\kwd{Sobol' Indices}
\kwd{Polynomial Chaos Approximation}
\end{keyword}


\end{frontmatter}

\tableofcontents

\section{Introduction}
\label{sec: Intro}

Computational models have penetrated everywhere. Over the past decades, they have become so complicated that there is an increasing need for special methods of their analysis. This analysis is even more important since the advent of artificial neural networks. In the context of AI safety, we have faced the interpretability problem \cite{Zhang2018}: how can we explain the decision made by the neural network? For instance, such methods of analysis are especially important in medical applications, where almost every decision is accompanied by the most serious consequences \cite{Begoli19}.

Being an important tool for investigation of computational models, \emph{sensitivity analysis} tries to find how different model input parameters influence the model output, what are the most influential parameters, and how to evaluate such effects quantitatively \cite{Saltelli2008}. Sensitivity analysis allows to understand the behavior of computational models better. Particularly, it allows us to separate all input parameters into {\it important (significant)}, {\it relatively important} and {\it unimportant (nonsignificant)} ones. Important parameters, {\it i.e.} parameters whose variability has a strong effect on the model output, may need to be controlled more accurately. As complex computational models often suffer from over-parametrization, by excluding unimportant parameters, we can potentially improve model quality, reduce parametrization and computational costs \cite{ElstatHastie2009}.


At its core, sensitivity analysis deals with the variation of the model response caused by the variability of model input parameters, where the latter may be controlled by the experimenter (adaptive design of experiments, active learning setting) or given {\it a-priori} (fixed sample \cite{Plischke13}). As an example, in \cite{Zeiler2014} image classification with Convolutional Neural Networks is considered, and by occluding different portions of the input image with a grey square and monitoring the output of CNN, it is revealed which parts of the scene are essential for the~classification.

Sensitivity analysis includes a wide range of metrics and techniques including the Morris method \cite{Morris1991}, linear regression-based methods \cite{Iooss2015}, variance-based methods \cite{Saltelli2008}, etc. See for details reviews  \cite{Iooss2015,Wei15}. Among all the metrics, we focus on \emph{Sobol' (sensitivity) index}, a common way to evaluate the nonlinear influence of model parameters \cite{Sobol93,Saltelli2010}. Sobol' indices quantify which portions of the output variance are explained by different input parameters and combinations thereof.

The approaches for the evaluation of Sobol' indices are usually divided into {\it Monte Carlo} and {\it metamodeling approaches}. Monte Carlo approaches run the analyzed model and conduct high-dimensional numerical integration to estimate Sobol' indices using direct formulas such as given in \cite{Tissot15}, SPF scheme \cite{SPFSobol2001}, formulas of Kucherenko \cite{Kucherenko11}, Myshetzskay \cite{Myshetskaya07}, Owen \cite{Owen12}, etc. They are relatively robust \cite{Yang2011}, but require a large number of model runs for an accurate estimation of each index; see \cite{DaVeiga13,Gamboa2016,Rugama18,Tissot15} for more information on the accuracy and convergence of these methods. Thus, Monte Carlo approaches are impractical for a number of applications, where each model evaluation has a high computational cost. 

On the other hand, {\it metamodeling approaches} allow one to reduce the required number of model runs \cite{Iooss2015}. Following this approach, we replace the original computational model with an approximating {\it metamodel} (also known as {\it surrogate model} or {\it response surface}) and use this approximation to calculate Sobol' indices. Unlike the original model, the metamodel is easy to simulate and has a known internal structure. Commonly used metamodels include Polynomial Chaos Approximation \cite{Sudret15}, Gaussian process (GP) approximation \cite{Marrel2009}, local polynomials \cite{DaVeiga06} and others.

\subsubsection*{Confidence of metamodeling approaches}

Following metamodeling approaches, we face the question of confidence in the results obtained: can we guarantee the closeness of the true Sobol' index and its estimate, especially for a small training sample? What size of the training sample is needed? Although the metamodeling approach can be more efficient, its accuracy is more difficult to analyze compared to Monte Carlo. Indeed, even though procedures like {\it cross-validation} \cite{ElstatHastie2009,Stone1974} allow to evaluate the quality of metamodel, the accuracy of complex statistics (such as Sobol' indices), derived from this metamodel, has a complicated dependence on the metamodel structure and its quality. In general, the approach includes two sources of uncertainty: {\it the metamodel error} and {\it the sampling (Monte Carlo) error}, where the latter is related to the calculation of Sobol' index from the metamodel (unless it is calculated analytically).

Several methods are proposed in the literature to solve this confidence problem, starting from a simple numerical simulation to more theoretical considerations. In \cite{CIbootstrapSobol}, bootstrap-based confidence intervals are constructed for Sobol' indices calculated from the Polynomial Chaos Approximation, and adaptive design is proposed for the computation of Sobol' indices with a given accuracy.  \cite{Marrel2009} proposes to approximate the original model with the Gaussian process metamodel and to use the metamodel simulation to estimate Sobol' indices. In addition, the probability distribution of the estimated Sobol' indices is obtained via numerical integration of  the Gaussian process realizations generated from the trained metamodel. According to another method, the sampling error is estimated with the bootstrap technique, and the metamodel error can be obtained as a solution of the multidimensional optimization problem, provided pointwise error bounds for the metamodel  are given \cite{Janon14UncAssesm}. An extension of this approach to the Gaussian process metamodel is presented in \cite{Gratiet14}.


Some authors consider the asymptotic approximation of Sobol' indices distribution and confidence intervals using the $\delta$-method \cite{OehlertDelta1992} in the classical regression setting of non-random design and the presence of random noise in the analyzed model output \cite{BurnaevPanin15,BurnaevPanin17,Xu2011}. \cite{Pronzato18} also gives an exact distribution of Sobol' indices estimates in this setting. A similar idea, including the asymptotic approximation through the $\delta$-method, is presented in \cite{Janon14AssympNorm}.  In addition, asymptotic properties of conditional moments estimation based on local polynomials are considered in~\cite{DaVeiga06}.

We should also mention \cite{Huang98}, where a general theory on rates of convergence of the least squares estimate is proposed, including the asymptotic convergence of ANOVA components. In addition, general adaptive estimators of quadratic functionals are considered in \cite{Laurent2000}.

\subsubsection*{Contribution}

Although the estimation of Sobol' indices with metamodels is widespread in practice, it lacks a theoretical justification. In this paper, we address the questions: what is the accuracy of metamodels-based Sobol' indices, the~risk of the estimate, and its rate of convergence? 
 
Following these questions, we establish a relation between the metamodel error and Sobol' indices errors; propose a method for indices quality control; and obtain the risk of Sobol' indices evaluation in the random design setting for a general class of metamodels. 

We assume that a target function is an element of some Hilbert space, and the approximation lies in its finite-dimensional subspace (the linear span of regressors). Two methods are exploited to estimate metamodels parameters: ordinary least squares and the projection (quasi-regression) method. In a non-parametric setting, we study the asymptotic and non-asymptotic behavior of Sobol' indices of all orders, including total-effects. The convergence rates and asymptotic selection of metamodel complexity are studied for different smoothness classes and for various metamodels that use Legendre, Chebyshev, and Trigonometric polynomials. 

In addition, we pay special attention to the case of noiseless observations in which the model input parameters completely determine the output of the model {\it i.e.} the variance $\V(y|\xb)=0$, and which often takes place in practice. 

Table~\ref{tab: SI_risk_all} summarizes the asymptotic results when the dimension of the approximation subspace grows with the increasing sample size. Our findings indicate that the absence of random noise in the output of the analyzed function, along with its high smoothness\footnote{See formal Definition~\ref{def: p-smoothness} of smoothness in Section~\ref{sec: Interpretation}.} and low dimension, is one of the main factors that provides the possibility of fast convergence rates of Sobol' indices estimates. Our results also give an opportunity to analyze various experimental designs used in global sensitivity analysis based on metamodels.

The paper is organized as follows: 
in Section \ref{sec: Sobol indices}, we review the definition of Sobol' indices, describe their evaluation using metamodels with tensor structure, and introduce two methods for the estimation of metamodels parameters. 
In Section \ref{sec: Risk}, the main results are presented, including a general relation between the accuracy of arbitrary metamodel and the error of estimated indices, and specific relations for the two methods. In addition, a method for quality control is proposed. In Section \ref{sec: Interpretation}, the obtained results are considered from an asymptotic point of view.  Experimental results\footnote{Python code is available at \url{http://www.github.com/Ivanx32/sobol_indices}} are provided in Section~\ref{sec: Experiments}. Proofs are given in Appendix~\ref{sec: Proofs}.


\setcellgapes{4pt}

\begin{table}[t]
\makegapedcells
\centering
\caption{\label{tab: SI_risk_all}Asymptotic bounds for the quadratic risk of Sobol' and total-effect indices estimates, depending on sample size $n$, input dimension $d$, and smoothness $p$.
}
\begin{tabular}{|c|c|c|c|c|}
\hline
\multicolumn{2}{|c|}{\textbf{}}        & \textbf{Legendre} & \textbf{Chebyshev} & \textbf{Trigonometric} \\ \hline
\multicolumn{2}{|c|}{measure}          & Uniform           & Arcsine            & Uniform                \\ \hline
\multirow{2}{*}{\rotatebox{90}{noiseless}} & LS   &   $\left(\frac{n}{\ln n}\right)^{-p/d}  $  &   $\left(\frac{n}{\ln n}\right)^{-2p/d}  $   &   $\left(\frac{n}{\ln n}\right)^{-2p/d}  $  \\[0pt] \cline{2-5} 
                          & projection & \multicolumn{3}{c|}{$n^{-\frac{2p}{2p+d}}$ }                                   \\[0pt] \hline
\multirow{2}{*}{\rotatebox{90}{noisy}}  & LS & \multicolumn{1}{c|}{\begin{tabular}[c]{@{}c@{}}$n^{-\frac{2p}{2p+d}}$, \;\; $p/d > 1/2$\\ $\left(\frac{n}{\ln n}\right)^{-p/d}$,  $\, p/d \leq 1/2$\end{tabular}}  &  $n^{-\frac{2p}{2p+d}}$    &  $n^{-\frac{2p}{2p+d}}$    \\ \cline{2-5} 
                          & projection &    \multicolumn{3}{c|}{$n^{-\frac{2p}{2p+d}}$ }    \\[0pt] \hline
\end{tabular}%
\end{table}


\section{Sobol' indices}
\label{sec: Sobol indices}

\subsection{Sobol-Hoeffding decomposition}
\label{subsec: Sobol-Hoeffding}

Consider a function $y = f(\xb)$, where $\xb = (x_1, \ldots, x_d) \in \XX \subseteq \R^d$~is a vector of \emph{input variables}, $y \in \R^1$ is an \emph{output variable}. \emph{Design space} $\XX$ is supposed to have a form $\XX = \XX_1 \times \ldots \times \XX_d$ with $\XX_i \subseteq \R$. 

Let us assume that there is a prescribed\footnote{The probability measure is fixed and known {\it a priori}.} probability measure  $\mu$ on the~design space, having the form of tensor product: 
$d \mu(\xb) = \otimes_{i=1}^d d\mu_i(x_i)$, where $\mu_i$ is a probability measure over $\XX_i$. The corresponding distribution represents the uncertainty and/or variability of the input variables, modelled as a random vector $\xb = (x_1, \ldots, x_d)$ with independent components. In this setting, the model output $y = f(\xb)$ becomes a stochastic variable.

Suppose that the function $f$ is square-integrable w.r.t. $\mu$ {\it i.e.} $f$ lies in Hilbert space $L^2(\XX, \mu)$ of real-valued functions\footnote{To be precise, $L^2(\XX, \mu)$ includes equivalence classes of functions that coincide $\mu$-almost everywhere.} on $\XX$ square-integrable for $\mu$. We have the following unique Sobol-Hoeffding decomposition of the model output \cite{Sobol93,Vaart}  given by
\begin{eqnarray*}
f(\xb) 
&=& 
f_0 + \sum_{i=1}^{d} f_i(x_i) + \sum_{1\leq i < j \leq d} f_{ij}(x_i, x_j) + \ldots + f_{1 \ldots d}(x_1, \ldots, x_d) 
\\
&=&  
\sum_{\ub \subseteq \{1, \ldots, d\}} f_\ub(\xb_\ub),
\end{eqnarray*}
with $2^d$ terms which satisfy
\begin{eqnarray*}
\E_i [f_\ub]  \eq \E_{\mu_i} [f_\ub] 
= 
\int_{\XX_i} f_{\ub}(\xb_\ub) d\mu_i(x_i) 
= 
0 \; \; \text{for} \;\; \forall i \in \ub,
\end{eqnarray*}
where $\ub \subseteq \{1, 2, \ldots, d\}$ is an index set, $\xb_{\ub}$ is the vector with components $x_i$ for $i \in \ub$, and $f_{\emptyset} \eq f_0 = \E_{\mu}[f(\xb)]$. As a consequence, for $\ub, \vb \subseteq \{1, 2, \ldots, d\}$
\[
\E_{\mu} \big[f_\ub(\xb_\ub)f_\vb(\xb_\vb) \big] 
= 
0 \;\; \text{if} \;\; \ub \neq \vb.
\]

Due to orthogonality of the summands, we can decompose the variance of the model output:
\begin{equation}
\label{eq: partial_var_hoeffding}
D = 
\V_{\mu}[f(\xb)] 
= 
\sum_{\substack{\ub \subseteq \{1, \ldots, d\}}  } \V_{\mu}[f_\ub(\xb_\ub)] 
= 
\sum_{\substack{\ub \subseteq \{1, \ldots, d\}}  } D_\ub.
\end{equation}
In this expansion, $D_\ub \triangleq \V_{\mu}[f_\ub(\xb_\ub)]$ is the contribution of the summand $f_\ub(\xb_\ub)$ to the output variance, also known as the \emph{partial variance}.

Define for $\V_{\mu}[f] > 0$ {\it Sobol' sensitivity index}  $S_{\ub}$.
\begin{definition}
\label{def: si}
    Sobol' index  of the subset $\xb_\ub,\; \ub \subseteq \{1, \ldots, d\}$ of model input variables is defined as
\begin{equation}
\label{eq: si_def_1}
	S_{\ub} = \frac{D_\ub}{D}.
\end{equation}
\end{definition}

\noindent Note that $S_{\ub} \in [0, 1]$ and $\sum_{\ub }S_{\ub} = 1$. 

Denote $d \mu_{\ub} \eq \otimes_{i \in \ub} d\mu_i$, $\E_{\ub} \eq \E_{\mu_{\ub}}$, $\V_{\ub} \eq \V_{\mu_{\ub}}$ and $\unsim\ub \eq  \{1, \ldots, d\} \backslash \ub$ (the subset of all input variables except the  variables belonging to $\ub$). Then for $\ub = \{i\}$ the  sensitivity index (\ref{eq: si_def_1}) can be calculated as
\[
   S_{i}  = \frac{\V_i \big[\E_{\unsim i} (f(\xb)| x_i)  \big]}{\V_{\mu}[f]},  \;\; i = 1, \ldots, d, \;\;
\]
and for $\ub = \{ij\}$ 
\[
S_{ij}   = \frac{\V_{ij} \big[\E_{\unsim ij} (f(\xb)| x_i, x_j)  \big]}{\V_{\mu}[f]} - S_i - S_j,  \;\;\;  i,\;j = 1, \ldots, d, \;\; i \neq j.
\]

We also define the quantity that characterizes the ``total'' contribution of variables $\xb_\ub$:  {\it total-effect} index also known as {\it total Sobol' index}.

\begin{definition}
     Total-effect index  of the subset $\xb_\ub,\; \ub \subseteq \{1, \ldots, d\}$ of model input variables is defined as
\begin{equation}
\label{eq: def_total_ind}
    T_{\ub} 
    = 
    \sum_{\ub \cap \vb \neq \emptyset} S_{\vb}   
    =
     1 - \sum_{\ub \cap \vb = \emptyset} S_{\vb} 
     =
    \frac{\E_{\unsim\ub} \big[\V_{\ub} (f(\xb)| \xb_{\unsim\ub})  \big]}{\V_{\mu}[f]}.
\end{equation}

\end{definition}

\noindent Note that $T_{\ub} \in [0, 1]$ and $ T_{\{1, \ldots, d\}} = \sum_{\ub}S_{\ub} = 1$. We can also formally define $S_{\emptyset} = T_{\emptyset} = 0$.


The question is how to efficiently calculate Sobol' indices?

\subsection{Metamodels-based sensitivity analysis}
\label{sec: metamodel}

\subsubsection{Sobol' indices and metamodels}

Direct calculation of Sobol' indices leads to computationally expensive multidimensional integration. To simplify this problem using the metamodeling approach, one can replace the original function $f(\xb)$ with the approximation $\hat{f}(\xb)$ that is better suited for computing of Sobol' indices. 

In general, metamodels-based sensitivity analysis includes the following steps: (a) the selection of  {\it the design of experiments}; (b) generation of responses of the analyzed model; (c) selection and construction of the metamodel based on the obtained training sample,  including its accuracy assessment; (d) the estimation of Sobol' indices using the constructed metamodel. The last step is based either on a known internal structure of the metamodel or on Monte Carlo simulation of the metamodel itself that is computationally cheap.

As such an approximation, we consider an expansion in a series of  $\mu$-orthonormal functions having a form of tensor product. Special cases of such approximation are Polynomial Chaos Approximation, low-rank tensor approximations (LRA) \cite{Konakli16}, Fourier approximation \cite{Popinski1999}, etc. Besides, by using Karhunen-Lo{\`e}ve decomposition one can reduce to this expansion Gaussian process approximation \cite{Pronzato18}. Similar generalized Chaos Expansions, built on general tensor Hilbert basis,  are considered in \cite{Roustant2020}. 

Let us formally introduce the approximation of this type. Denote scalar product and norm for $g, h \in L^2(\XX, \mu)$ as 
$\la g, h \ra_{\mu} = \int_{\xb \in \XX} g(\xb)h(\xb)d\mu(\xb) $,  $\|g\|^2_{\mu} \triangleq \|g\|^2_{L^2(\XX, \mu)}=\int_{\xb \in \XX} g^2(\xb)d\mu(\xb)$. In addition, {\it the uniform norm} $\|g\|_{L^{\infty}} \eq \|g\|_{L^{\infty}(\XX)} \eq \sup_{\xb \in \XX} |g(\xb)|$. The norm in Euclidean vector spaces is denoted by $\|\cdot\|$.

Define the regressors. Suppose there exists a function set $\{\Psi_{\al}(\xb) \}$ in $L^2(\XX, \mu)$ parametrized with multi-index\footnote{We use the notation $\NN \eq \{0,1,2,\ldots\}$ for the set of all nonnegative integers and $\NN_{+}$ for positive integers.} $\al = (\alpha_1, \ldots ,\alpha_d) \in \NN^d$ that consists of $\mu$-orthonormal functions having a form of tensor product of $d$ $\mu_i$-orthonormal univariate function families $\{\psi_{\alpha_i}^{(i)}, \; \alpha_i \in \NN\}$ with $\psi_{0}^{(i)} \eq 1$ and $\ee_{i}\big[\psi^{(i)}_{\alpha}(x_i) \psi^{(i)}_{\beta}(x_i) \big] 
= 
\delta_{\alpha\beta} \; \text{for}\; \alpha, \beta \in \NN$, where $\delta$ is the Kronecker symbol. As~a~result,
\begin{equation}
\label{eq: ONB}
\begin{aligned}
\Psi_{\mathbf{0}}(\xb) \eq 1, \;\;
\Psi_{\al}(\xb) 
&= 
\prod_{i=1}^d \psi_{\alpha_i}^{(i)}(x_i),
\;\;
\xb \in \XX,  
  \\
\la \Psi_{\al},\Psi_{\be}\ra_{\mu}
&= 
\delta_{\al\be}, \;\; \al, \be \in \NN^d. 
\end{aligned}
\end{equation}

An example of such set for a uniform distribution is multivariate Legendre polynomials which are additionally normalized to have a unit variance w.r.t. a uniform measure on $[-1,1]^d$ (see Section~\ref{sec: legendre_asymp_int}). 

\begin{remark}
\label{rem: finite_set_discrete_measure}
W.l.o.g. we will consider the set $\{\Psi_{\al}(\xb)\}$ that  consists of an infinite number of elements, but all further statements also remain  true for a finite number. For example, such a case can take place for the discrete measure  $\mu$. We only assume that there is at least one non-constant element $\psi_{1}^{(i)}$ for each input dimension $i = 1, \ldots, d$ that leads to $\{\Psi_{\al}(\xb)\}$ which contains at least the following $2^d$ elements: $\Psi_{(0, 0, \ldots, 0)}(\xb) =1$, $\Psi_{(1, 0, \ldots, 0)}(\xb) = \psi_{1}^{(1)}(x_1)$,  $\Psi_{(0, 1, \ldots, 0)}(\xb) = \psi_{1}^{(2)}(x_2)$, \ldots, $\Psi_{(1, 1, \ldots, 1)}(\xb) = \prod_{i=1}^d \psi_{1}^{(i)}(x_i)$.
\end{remark}

Define the metamodel with tensor structure. Let the approximation be the linear combination of $N$ functions from the set $\{\Psi_{\al}(\xb), \; \al \in \rr_N\}$ for some $\rr_N \subset \NN^d$. In other words, 
\begin{equation}
\label{eq: hat_f_expansion}
\hat{f}(\xb) = \sum_{\al \in \rr_N} \hat{c}_{\al} \Psi_{\al}(\xb), \;\;\; \xb \in \XX, \;\; \hat{c}_{\al} \in \R.
\end{equation} 

If the regressors $\Psi_{\al}$ are $\mu$-orthogonal polynomials, then (\ref{eq: hat_f_expansion}) corresponds to Polynomial Chaos Approximation.

One of the important advantages of the presented approximation type is that it allows to calculate Sobol' indices analytically from the expansion coefficients. Indeed, it can be shown (see \cite{Pronzato18}) that for the function $\hat{f}$ defined as (\ref{eq: hat_f_expansion}) using the  orthonormal set (\ref{eq: ONB}) it holds
\begin{equation}
\label{eq:norm_var_hat_f}
\|\hat{f}\|^2_{\mu} = \sum_{\al \in \rr_N} \hat{c}_{\al}^2, \;\; \E_{\mu} [\hat{f}] =  \hat{c}_{\mathbf{0}}, \;\; \V_{\mu} [\hat{f}] =  \sum_{\al \in \rr_N\backslash\mathbf{0}} \hat{c}_{\al}^2,
\end{equation}
and if $\V_{\mu}[\hat{f}] > 0$, then Sobol' index for $\xb_\ub$ variables of  $\hat{f}(\xb)$ reads
\begin{equation}
\label{eq:sobol_tensor}
\hat{S}_\ub(\hat{\cb}) = \frac{ \sum_{\al \in \T_\ub} \hat{c}_{\al}^2}{\sum_{\al \in \rr_N\backslash\mathbf{0}} \hat{c}_{\al}^2 }, \;\; \ub \subseteq \{1, \ldots, d\}, \; \ub \neq \emptyset,
\end{equation}
where $\T_\ub \eq \T_\ub[\rr_N]$ is the subset of  $\rr_N$ that consists of such multi-indices that only indices corresponding to variables $\xb_\ub$ are nonzero: $\T_\ub =  \{\al \in \rr_N\colon \; \alpha_i > 0 \; \text{for all} \; i \in \ub; \;  \alpha_i = 0  \; \text{for} \; i \notin \ub \}$. E.g.~$\T_i[\rr_N] = \big \{\al \in \rr_N\colon \al = (0, \ldots, \alpha_i, \ldots, 0)$, $\alpha_i>0 \big\}$.  The vector of coefficients $\{\hat{c}_{\al}, \; \al \in \rr_N\}$ is denoted as $\hat{\cb} \in \R^N$.

Similarly, the total-effect index reads
\begin{equation}
\label{eq:total_sobol_tensor}
\hat{T}_\ub(\hat{\cb}) = \frac{ \sum_{\al \in \TT_\ub} \hat{c}_{\al}^2}{\sum_{\al \in \rr_N\backslash\mathbf{0}} \hat{c}_{\al}^2 }, \;\; \ub \subseteq \{1, \ldots, d\}, \; \ub \neq \emptyset,
\end{equation}
where $\TT_\ub \eq \TT_\ub[\rr_N]$ is the subset of  $\rr_N$ that consists of such multi-indices that at least one index corresponding to variables $\xb_\ub$ is nonzero: 
$\TT_\ub= \{\al \in \rr_N\colon \; \sum_{i \in \ub}\alpha_i > 0\}$.

\begin{remark}
Similar analytical expressions of Sobol' indices for Polynomial Chaos was obtained in \cite{Sudret}, and for Fourier approximation such an approach is exploited in {\it Fourier amplitude sensitivity test (FAST)} \cite{FASTCukier1978}  and {\it Random balance designs (RBD)} \cite{Tarantola06}.  
\end{remark}

In the next section, we consider the strategies for constructing an approximation based on a finite training sample.


\subsubsection{Approximation construction}
\label{sec: Appr_constr}

Consider some fixed nested sets of $d$-dimensional  multi-indices 
\begin{equation}
\label{eq: fixed_multiind_order}
\{\mathbf{0}\} = \rr_1 \subset \ldots \subset \rr_N \ldots \subset \rr_{\infty}= \NN^d,
\end{equation}
where $|\rr_N| = N$ for all $N \in \NN_{+}$. We will refer to every $\rr_N$ as {\it the truncation set}, as we assume that all expansion coefficients of the corresponding approximation, not belonging to $\rr_N$, are zero: $\hat{c}_{\al} \eq 0$ for $\al \notin \rr_N$. Selecting $N$, we define in some sense the complexity of the approximating model.

Denote the subspace of all linear combinations of $\{\Psi_{\al}, \; \al \in \rr_N\}$ as $\pnu_N \eq span\{\Psi_{\al}, \; \al \in \rr_N\}$.  {\it Theoretical orthogonal projection} of $f$ onto $\pnu_N$ w.r.t. $\mu$-norm is defined as
\begin{equation}
\label{eq: best_apprx_pnu_N}
f_N \eq \argmin_{\hat{f} \in \pnu_N}\|f - \hat{f}\|_{\mu}.
\end{equation} 
Define also {\it the residual}
\begin{equation}
\label{eq: e_N}
e_N(\xb) \eq f(\xb) -  f_N(\xb),
\end{equation}
then the error of the best approximation of $f$ in the space $\pnu_N$ w.r.t. $\mu$-norm is $\|e_N\|_{\mu}$. We have the following representation of arbitrary $f \in L^2(\XX, \mu)$:
\begin{equation}
\label{eq: f_full_representation}
f(\xb) = e_{\infty}(\xb) + \sum_{\al \in \NN^d} c_{\al} \Psi_{\al}(\xb), \;\; \text{with}\;\;\sum_{\al \in \NN^d} c_{\al}^2 < \infty.
\end{equation}




{\bf The observation model.} In general, the only information about $f$  comes from the observations. We suppose that for some \emph{design of experiments} $ \des = \{\xb_i \in \XX \}_{i = 1}^n \in  \R^{n \times d}$ in  {\it the samples space} $\XX^n$ we can obtain a set of model responses and form \emph{a~training sample}:
\begin{equation}
\label{eq:training_sample}
\mathcal{S} = \big(\xb_i, \; y_i = f(\xb_i) + \eta_i \big)_{i = 1}^n, \;\;
\end{equation}
where $\eta_i$ are i.i.d. realizations of random noise $\eta$ with $\E \eta = 0$ and $\V \eta = \sigma^2 < \infty$, independent from $\xb$. In the matrix form
\begin{equation*}
\label{eq:training_sample_matrix}
\mathcal{S} = \big(\des \in \XX^n, \;  Y = f(\des) + \etab \in \R^n \big), 
\end{equation*}
where $f(\des) \eq \big (f(\xb_1), \ldots, f(\xb_n) \big)^T$ and $\etab = \big (\eta_1, \ldots, \eta_n)^T $. We will consider both the general case of {\it noisy} observations and the special {\it noiseless} case corresponding $\sigma^2 = 0$. Noiseless case may also include situations when we can control the state of random generator inside the analyzed model.

Thus, we cannot obtain projection $f_N$ directly but can try to do this numerically. We will describe two basic ways to estimate the expansion coefficients in~(\ref{eq: hat_f_expansion}) from the training sample $\mathcal{S}$ of size $n$.

{\bf Projection{\footnote{This name may be misleading but allows us to distinguish between the methods.} method}} is based on quasi-regression \cite{Jian2001} and is used in the RBD-like methods \cite{Tarantola06}. Starting from the representation (\ref{eq: f_full_representation}) of $f$, we have
\[
c_{\al} = \la f,  \Psi_{\al}\ra_{\mu},
\;\;
\al \in \NN^d.
\]

Assuming the design is randomly i.i.d. sampled from the probability measure $\mu$, one can estimate the integral $\la f,  \Psi_{\al}\ra_{\mu}$ numerically. Replacing the scalar product with its empirical approximation based on the training sample, obtain
\begin{equation}
\label{eq: Projection_coeff}
\hat{c}_{\al}  
=
\frac{1}{n} \sum_{i=1}^n y_i \Psi_{\al}(\xb_i)
=
\frac{1}{n}  Y^T \Psi_{\al}(\des), \;\; \al \in \rr_N,
\end{equation}
where $\Psi_{\al}(\des) \eq \big(\Psi_{\al}(\xb_1), \ldots,  \Psi_{\al}(\xb_n) \big)^T \in \R^N$. We have
\[
\hat{\cb} 
= 
\frac{1}{n} \, \Phi^T Y \in \R^{N},
\]
with {\it the design matrix} $\Phi = \boldsymbol \Psi(\des) = \big\{\Psi_{\al}(\des), \; \al \in \rr_N \big\} \in \R^{n \times N}$.

{\bf Ordinary least squares (OLS)} method is often used for Polynomial Chaos Approximation construction \cite{Sudret15} and in many other cases. Replace the theoretical norm in (\ref{eq: best_apprx_pnu_N}) with its empirical version:
\begin{equation}
\label{eq: ls_empirical_minimization}
\hat{f}^{LS} = \argmin_{\hat{f} \in \pnu_N} \sum_{i=1}^n \big[y_i - \hat{f}(\xb_i) \big]^2.
\end{equation}
Assuming $\det (\Phi^T \Phi) \neq 0$, this minimization problem leads to
\begin{equation}
\label{eq: ls_coeff}
\hat{\cb} = (\Phi^T \Phi)^{-1}\Phi^T Y.
\end{equation}

We will refer to the approximations constructed based on these  two methods as $\hat{f}^{P}$ and  $\hat{f}^{LS}$ correspondingly. Related Sobol' indices, estimated via these two  approximations using (\ref{eq:sobol_tensor}), are denoted as $\hat{S}^{P}$ and  $\hat{S}^{LS}$ correspondingly. The next question is how the quality of approximation limits the quality of indices.

\section{Risk of Sobol' indices estimation}
\label{sec: Risk}

%

We begin with a general outline of this section. At first, we consider the relationship of the distance between the function and its arbitrary approximation and the distance between corresponding Sobol' indices (using various distance metrics). Based on this relation, we propose a method for quality control of indices estimates. In other words, the idea of the method is to establish a connection between the error of estimated Sobol' indices that is difficult to measure, and the relatively easily measurable approximation error. In addition, we show that the presented error bounds are achievable and, in this sense, cannot be improved.

The second part of the section is devoted to the risk of Sobol' indices estimates in the random design setting.  We generalize the obtained results for the approximations having tensor structure, which were constructed based on the projection and least squares methods with randomly chosen design points.

All further results are valid for both Sobol' indices and total-effects of all orders unless otherwise stated. In order to avoid duplication, we use the notation $S_\ub$ for indices of both types in theorems' statements.

\subsection{Error bounds}
\label{sec: error_bounds}

\subsubsection{Main relationships}

Our first goal is to assure that the closeness in some sense of the function and its arbitrary approximation $\hat{f} \approx f$ leads to the closeness of their (total) Sobol' indices $\hat{S}_{\ub}$ and $S_{\ub}$.  Note that the opposite is not true in general. To characterize the closeness of functions, we will use the relative error of approximation w.r.t. $\mu$-norm given by 
\begin{equation}
\label{eq: def_relative_eps}
\mer \eq \frac{\|f - \hat{f}\|_{\mu} }{ \V^{1/2}_{\mu}[f]}. 
\end{equation}

\begin{theorem}
\label{th: si_error_general}
For any $f, \hat{f} \in L^2(\XX, \mu)$ such that $\V_\mu[f] > 0$, $\V_\mu[\hat{f}] > 0$, it holds for corresponding Sobol' and total-effect indices for $\ub \subseteq \{1, \ldots, d\}$
\begin{eqnarray}
\label{eq: si_error_general_each_S}
\big|S_\ub - \hat{S}_\ub \big|
&\leq&  
\left\{ \sqrt{S_\ub(1-\hat{S}_\ub)} + \sqrt{\hat{S}_\ub(1-S_\ub)}\right\} 
\cdot 
\mer,
\\
\label{eq: si_error_general}
\max_{\ub}\big|S_\ub - \hat{S}_\ub \big|
&\leq&  
\mer.
\end{eqnarray}

\end{theorem}



The right-hand side of (\ref{eq: si_error_general_each_S}) contains both $S_\ub$ and $\hat{S}_\ub$. Solving this inequality w.r.t. $\big|S_\ub - \hat{S}_\ub \big|$ and simplifying the solution, obtain  an error estimate for the Sobol' index that depends only on the relative error $\mer$ and the true index $S_\ub$.
\begin{corollary}
\label{th: si_error_general_resolved}
For any $f, \hat{f} \in L^2(\XX, \mu)$ such that $\V_\mu[f] > 0$, $\V_\mu[\hat{f}] > 0$, it holds for corresponding Sobol' and total-effect indices for $\ub \subseteq \{1, \ldots, d\}$

\begin{equation}
\label{eq: si_error_general_each_S_approx}
\big|S_\ub - \hat{S}_\ub \big|
\leq  
\min\left(1, \; \mer + 2\sqrt{S_{\ub}}, \;  \mer + 2\sqrt{1-S_{\ub}} \right) \cdot \mer.
\end{equation}
In particular, assume that Sobol' or total-effect index $S_\ub \in \{0, 1\}$ for some $\ub$, then
\begin{equation}
\label{eq: si_error_general_zero-one}
\big|S_\ub - \hat{S}_\ub \big|
\leq  
\mer^2.
\end{equation}

\end{corollary}

As we see, the error bound of Sobol' indices is proportional to the fraction of the standard deviation not explained by the approximation, and one can expect a decrease in this error, respectively, with an increase in the quality of the approximation. 

The following corollary gives the bound for the sum of errors of Sobol' indices for all $2^d$ different subgroups of  variables (not valid for total-effects).
\begin{corollary}
\label{th: si_error_general_max}
For any $f, \hat{f} \in L^2(\XX, \mu)$ such that $\V_\mu[f] > 0$, $\V_\mu[\hat{f}] > 0$, it holds for corresponding Sobol' indices for $\ub \subseteq \{1, \ldots, d\}$
\begin{eqnarray}
\label{eq: si_error_general_sum_abs}
\sum_{\ub} \big|S_\ub - \hat{S}_\ub \big|
&\leq&  
2 \cdot \mer,
\\
\label{eq: si_error_general_sum_squared}
\sum_{\ub} \big(S_\ub - \hat{S}_\ub \big)^2
&\leq&
2 \cdot \mer^2.
\end{eqnarray}

\end{corollary}




\begin{example}
Let $\|f(\xb) - \hat{f}(\xb)\|_{\mu} \leq 0.1 \cdot \V^{1/2}_{\mu}[f]$ (for example, this is true if $\big|f(\xb) - \hat{f}(\xb)\big| \leq 0.1 \cdot \V^{1/2}_{\mu}[f]$ for $\xb \in \XX$). Then the errors of Sobol' and total-effect indices are bounded: $\max_{\ub} \{|S_{\ub} - \hat{S}_{\ub}|, |T_{\ub} - \hat{T}_{\ub}|\} \leq 0.1$. Besides, for Sobol' indices $\sum_{\ub}|S_\ub - \hat{S}_\ub | \leq  0.2$ and  $\sum_{\ub}(S_\ub - \hat{S}_\ub)^2
\leq  0.02$. Moreover, if for some variable groups $\ub, \vb \subseteq \{1, \ldots, d\}$ it holds $S_{\ub} = 0$ and $T_{\vb} = 0$, then $\max \{\hat{S}_{\ub}, \hat{T}_{\vb}\} \leq 0.01$.
\end{example}

Obviously of course, the closeness of all Sobol' indices of two functions does not necessarily lead to the closeness of functions. For example, consider a uniform measure $\mu$ on $[0,1]^2$ and two functions: $f(x_1, x_2) = x_1$ and $\hat{f}(x_1, x_2) = \sin 2\pi x_1$. Although all Sobol' indices (including total) of these functions  are the same, the functions are not equal in the sense of $\mu$-norm, {\it i.e.} $\|f - \hat{f}\|_{\mu} > 0$.


{\bf Discussion.} To the best of our knowledge, for the first time, the explicit relation of the errors of Sobol' indices (total-effects) and the quality of the approximation was obtained. Our results lead directly to a simple, practical method for the estimation of indices errors (see section \ref{sec: a_new_method}).

For comparison, we can consider the approach described in \cite{Janon14UncAssesm}.  The method involves the construction of a metamodel $\hat{f}$ and the use of estimation of Sobol' indices $\hat{S}_\ub^{MC}[\hat f]$ by the Monte Carlo method\footnote{With a large computational budget compared to $n$.}. In order to assess the metamodel error of the index, one can find bounds for the Monte Carlo estimator applied to the original model, {\it i.e.} to find for each $\ub \subseteq \{1, \ldots, d\}$ the values $\hat S_{\ub}^{m}$ and $\hat S_{\ub}^{M}$ (both depending on~$\hat{f}$) such that 
\[
\hat S_{\ub}^{m} \leq \hat{S}_{\ub}^{MC}[f] \leq \hat S_{\ub}^{M}.
\]
The values $\hat S_{\ub}^{m}$ and $\hat S_{\ub}^{M}$ are the solutions of some optimization problems\footnote{For more details see \cite{Janon14UncAssesm}.} connected with the variability of the Monte Carlo estimator $\hat{S}_\ub^{MC}[f]$ given {\it the pointwise error bound} (the function $\eps(\xb)$ such that $\big|f(\xb) - \hat{f}({\xb})\big| \leq \eps(\xb)$ for each $\xb \in \XX$). Of course,  only specific metamodel types can provide such pointwise error; for example, Gaussian process approximation\footnote{To be precise, GP provides only probabilistic pointwise errors.}. In addition, to take into account the sampling error of the index, the bootstrap method is used, and the repetition of the optimization procedure for each bootstrap subsample is required. As a result, we have a computationally demanding scheme with no closed-form solution that relies on the accuracy of the pointwise error estimate.




In contrast to the previous approach, our method takes into account only the ``average'' accuracy of the approximation which is easier to estimate and which is more robust to outliers compared to pointwise accuracy. In addition, the presented bounds are not associated with the specific type of approximation (see also Remark~\ref{re: analytical_or_large_budget}). Thus, these results are valid for arbitrary approximating models such as Polynomial Chaos, Gaussian processes, local polynomials, and others. The results obtained can be especially helpful for methods such as FAST and RBD that use Fourier approximation, but not explicitly. Previously, there were no reliable estimates of the accuracy for these methods.

It should be noted that the obtained error bounds also do not depend on the specific design of experiments for the approximation construction, so they are also valid if the design is associated with a measure other than $\mu$ (or even in the case of adaptive design \cite{BurnaevPanin17,Pronzato18}). Following this line of thought, one can show that the metamodeling approaches can provide accurate estimates of Sobol' indices for various sampling strategies, which is especially important if we cannot control  input sampling. Indeed, suppose that the experimental design is sampled from an unknown probability measure $\nu$ on $\XX$ such that $d(\mu, \nu) \eq \sup_{f} \{\|f\|_{\mu}/\|f\|_{\nu}\}$ is finite (the supremum is taken over all functions square-integrable with respect to both $\mu$ and $\nu$). Based on~(\ref{eq: si_error_general}),
\begin{equation}
\big|S_\ub - \hat{S}_\ub \big|
\leq
\frac{\|f - \hat{f}\|_{\mu} }{ \V^{1/2}_{\mu}[f]}
\leq
\frac{d(\mu, \nu)}{\V^{1/2}_{\mu}[f]} \cdot \|f - \hat{f}\|_{\nu},
\end{equation}
and hence a sufficient decrease in $\|f - \hat{f}\|_{\nu}$ results in decreasing index error.

This setting is connected with the so-called {\it second-level global sensitivity analysis}, which aims to quantify the impact of the uncertainty of the input probability distribution \cite{meynaoui2019new}. For example, in \cite{Hart2019} the robustness of the Sobol' indices to marginal distribution uncertainty is studied by calculating a derivative of the Sobol' index with respect to the probability density functions to determine an optimal perturbation which gives the largest change in the Sobol' index. The second-level analysis requires a very large number of model evaluations, and metamodeling approaches can help circumvent this problem. In this case, the obtained results provide the requirements for the approximation accuracy to ensure the closeness of Sobol' indices and their estimates even when the measure is changed:
\begin{equation}
\Big|S_\ub^{[\nu]} - \hat{S}_\ub^{[\nu]} \Big|
\leq
\frac{\|f - \hat{f}\|_{\nu} }{ \V^{1/2}_{\nu}[\hat f]}
\leq
\frac{d(\nu, \mu)}{\V^{1/2}_{\nu}[\hat f]} \cdot \|f - \hat{f}\|_{\mu},
\end{equation}
where $S_\ub^{[\nu]}$ and $\hat{S}_\ub^{[\nu]}$ are Sobol' indices of $f$ and $\hat{f}$ w.r.t the probability measure~$\nu$; $d(\nu, \mu)$ is supposed to be finite.

\subsubsection{A new method for quality control}
\label{sec: a_new_method}

Corollary~\ref{th: si_error_general_resolved} allows us to propose a new method for the quality control for Sobol' and total-effect indices assessment. From a practical point of view, it may be useful to represent the inequality for errors of indices (\ref{eq: si_error_general_each_S_approx}) in the form
\begin{equation}
\label{eq: si_error_general_each_S_practice}
\big|S_\ub - \hat{S}_\ub \big|
\leq  
\min\left(1, \; \mer_2 + 2\sqrt{\hat{S}_{\ub}}, \;  \mer_2 + 2\sqrt{1-\hat{S}_{\ub}} \right) \cdot \mer_2,
\end{equation}
with $\mer_2 \eq \|f - \hat{f}\|_{\mu} \, \cdot \, \min\big\{\V^{-1/2}_{\mu}[f], \V^{-1/2}_{\mu}[\hat f]\big\}$, which follows from the symmetry of Theorem~\ref{th: si_error_general} (and consequently Corollary~\ref{th: si_error_general_resolved}) w.r.t. $f$ and $\hat f$. 


Empirically, in practice to estimate $\mer_2$ one can replace the theoretical approximation error $\|f - \hat{f}\|_{\mu}$ and the variances $\V_{\mu}[f], \, \V_{\mu}[\hat f]$ with their sample estimates. Particularly, to replace the  true approximation error with {\it Root Mean Square Error (RMSE)}, obtained with procedures like {\it holdout method}. In the noiseless case ($\sigma^2 =0$) one can use the following estimate for the true error $\|f - \hat{f}\|_{\mu}$:
\begin{equation}
\label{eq: MSE_noiseless_app}
\text{RMSE}
=
\Bigg\{ \frac{1}{n_{test}} 
\sum_{i=1}^{n_{test}} \big[f(\xb_i) - \hat{f}(\xb_i) \big]^2 
\Bigg\}^{1/2}, 
\end{equation} 
where $\{\xb_i\}_{i=1}^{n_{test}}$ is a sample of size $n_{test}$ {\it i.i.d.} generated from the measure $\mu$ that is independent from the design $\des$ and is not used for the~construction of approximation $\hat{f}$. Note that some types of approximations $\hat{f}$ provide analytical expressions for their variances, see for example  (\ref{eq:norm_var_hat_f}).

\begin{remark} 
\label{re: analytical_or_large_budget}
One assumption needs to be emphasized. We suppose that the used metamodel allows the Sobol' indices to be calculated analytically, or at least that the computational budget is sufficient to provide a negligibly small error of this calculation. If this is not the case, one can adjust the presented bounds, for example, using the bootstrap method~\cite{Janon14UncAssesm}. 
\end{remark}

\begin{remark}
The obtained results can also be used as a justification for a thresholding rule. For example, according to (\ref{eq: si_error_general_zero-one}), $\hat{T}_{\ub} > \mer^2$ for some $\ub$ implies $T_{\ub} > 0$, and one can select a sample estimate of $\mer^2$ as the threshold that separates relatively important variables and their groups from completely unimportant. 
\end{remark}

\subsubsection{Tightness of error bounds}
\label{sec: Tightness_of_bounds}

We present an illustration of Theorem~\ref{th: si_error_general} and Corollary~\ref{th: si_error_general_max} and show that the corresponding bounds are achievable. Consider two functions on $\XX \subseteq \R^d$ having simple additive forms, constructed as linear combinations of two elements\footnote{See also Remark~\ref{rem: finite_set_discrete_measure}.} of the orthonormal set (\ref{eq: ONB}):
\begin{equation}
\label{eq: two_func_example_tightness}
\begin{aligned}
f(\xb) 
&= 
c_1 \cdot \Psi_{(1,0, \dots)}(\xb) +  c_2 \cdot \Psi_{(0,1, \dots)}(\xb), 
\\
\hat{f}(\xb) 
&= 
\hat{c}_1 \cdot \Psi_{(1,0, \dots)}(\xb) +  \hat{c}_2 \cdot \Psi_{(0,1, \dots)}(\xb),
\end{aligned}
\end{equation}
where coefficients $c_1, c_2, \hat{c}_1, \hat{c}_2 \geq 0$, $c_1^2 +c_2^2 > 0$, $\hat{c}_1^2 +\hat{c}_2^2 > 0$. Taking into account  $\Psi_{(1, 0, \ldots, 0)}(\xb) = \psi_{1}^{(1)}(x_1)$ and $\Psi_{(0, 1, \ldots, 0)}(\xb) = \psi_{1}^{(2)}(x_2)$, the corresponding Sobol' and total-effect indices for $f$ and $\hat f$:
\begin{eqnarray*}
S_1 &=& T_1 =  \frac{c_1^2}{c_1^2 + c_2^2}, \;\;  
S_2 = T_2 =  \frac{c_2^2}{c_1^2 + c_2^2}, \;\;
S_{\ub} = 0 \; \text{for}\; \ub \notin \big\{\{1\}, \{2\}\big\},
\\
\hat{S}_1 &=& \hat{T}_1 =  \frac{\hat{c}_1^2}{\hat{c}_1^2 + \hat{c}_2^2}, \;\;  
\hat{S}_2 = \hat{T}_2 =  \frac{\hat{c}_2^2}{\hat{c}_1^2 + \hat{c}_2^2}, \;\;
\hat{S}_{\ub} = 0 \; \text{for}\; \ub \notin \big\{\{1\}, \{2\}\big\}.
\end{eqnarray*}

\begin{figure}[t!]
\begin{center}
  \includegraphics[scale=1]{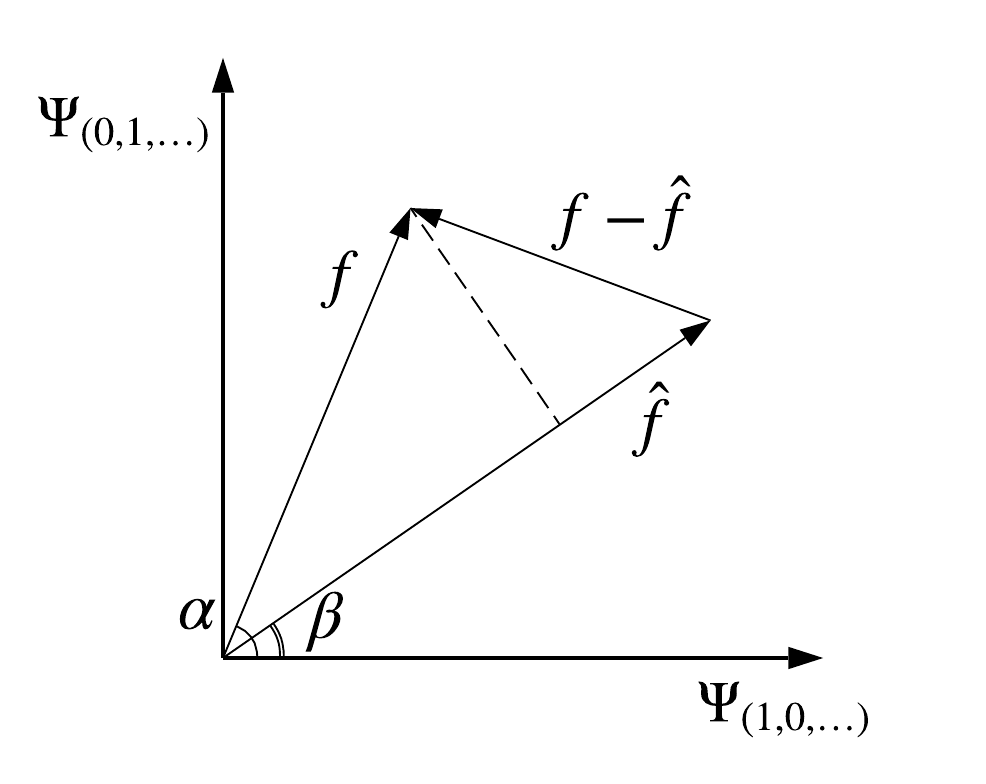}
  \caption{2D example for the functions (\ref{eq: two_func_example_tightness}), where $S_1 = \cos^2 \alpha$, $\hat{S}_1 = \cos^2 \beta$, and $\|f - \hat f\|_{\mu} \geq \|f\|_{\mu} \cdot |\sin(\alpha - \beta)|$.}
  \label{fig:2d_example}
\end{center}
\end{figure}

Compare $|S_1 - \hat{S}_1|$ and $\|f - \hat{f}\|_{\mu} = \sqrt{(c_1-\hat{c}_1)^2 + (c_2-\hat{c}_2)^2}$. Denote $\alpha \eq \arccos S_{1}^{1/2}$ and  $\beta \eq \arccos \hat{S}_{1}^{1/2}$ (see Fig.~\ref{fig:2d_example}). Notice that
\begin{eqnarray*}
\Big|S_1 - \hat{S}_1 \Big| 
&=&
\big|\cos^2 \alpha - \cos^2 \beta \big|
=
\sin(\alpha+\beta) \cdot \big|\sin(\alpha - \beta) \big|
\\
&=&
\left\{ \sqrt{S_1(1-\hat{S}_1)} + \sqrt{\hat{S}_1(1-S_1)}\right\} 
\cdot
\big|\sin(\alpha - \beta) \big|
\\
&\leq&
\left\{ \sqrt{S_1(1-\hat{S}_1)} + \sqrt{\hat{S}_1(1-S_1)}\right\} 
\cdot
\mer.
\end{eqnarray*}
Moreover, the last inequality turns into equality if we additionally require $\la f, \hat{f}\ra_{\mu} \neq 0$ and set $\la \hat{f}, f - \hat{f} \ra_{\mu} = 0$, see Fig.~\ref{fig:2d_example}. Thus, we can achieve the bound (\ref{eq: si_error_general_each_S}) for any $S_1, \hat S_1 \in [0, 1]$ except two cases corresponding to $\la f, \hat{f}\ra_{\mu} = 0$: $S_1=0$, $\hat S_1=1$ and $S_1=1$, $\hat S_1=0$. It can be shown that in these two cases the index error  
does not reach the upper bound  (\ref{eq: si_error_general_each_S}), but can be arbitrarily close to it.

One can summarize these arguments in the following theorem.
\begin{theorem}
\label{th: si_error_general_strict}
For any subset $\ub \subseteq \{1, \ldots, d\}$ and any values $S_{\ub}, \hat{S}_{\ub} \in (0, 1)$ there exists $f, \hat{f} \in L^2(\XX, \mu)$ having Sobol' (total-effect) indices for $\xb_{\ub}$~variables equal $S_{\ub}$ and $\hat{S}_{\ub}$ correspondingly, such that
\begin{eqnarray}
\label{eq: si_error_general_each_S_strict}
\big|S_\ub - \hat{S}_\ub \big|
&=&  
\left\{ \sqrt{S_\ub(1-\hat{S}_\ub)} + \sqrt{\hat{S}_\ub(1-S_\ub)}\right\} 
\cdot 
\mer.
\end{eqnarray}

\end{theorem}


Select the coefficients $c_1, c_2, \hat{c}_1, \hat{c}_2$ such that $\sin(\alpha+\beta) = 1$ and $\la \hat{f}, f - \hat{f} \ra_{\mu} = 0$. For example, we may take $\alpha=\pi/3$, $\beta=\pi/6$ and
\begin{equation}
\label{eq: example_strict_bound_degree30}
\begin{aligned}
f(\xb) 
&= 
\frac{1}{2} \cdot \Psi_{(1,0, \dots)}(\xb) +  \frac{\sqrt{3}}{2} \cdot \Psi_{(0,1, \dots)}(\xb),
\\
\hat{f}(\xb) 
&= 
\frac{3}{4} \cdot \Psi_{(1,0, \dots)}(\xb) +  \frac{\sqrt{3}}{4} \cdot \Psi_{(0,1, \dots)}(\xb).
\end{aligned}
\end{equation}

One can see that  $S_1 = \cos^2 \frac{\pi}{3} = 1/4$, $S_2  = 3/4$, $\hat{S}_1 = \cos^2 \frac{\pi}{6}= 3/4$, $\hat{S}_2 = 1/4$ and the bound (\ref{eq: si_error_general})  in Theorem~\ref{th: si_error_general} is achieved:
\begin{equation}
\label{eq: example_bound_max_3}
\max_{\ub}\left \{\big|S_\ub - \hat{S}_\ub \big|, |T_\ub - \hat{T}_\ub \big|\right\}
=
\big|S_1 - \hat{S}_1 \big|
= 
\mer.
\end{equation}
The bounds (\ref{eq: si_error_general_sum_abs}, \ref{eq: si_error_general_sum_squared}) in Corollary~\ref{th: si_error_general_max} are also achieved:
\begin{eqnarray*}
\sum_{\ub} \big| S_{\ub} - \hat{S}_{\ub} \big|
&=&
\big| S_{1} - \hat{S}_{1} \big|
+
\big| S_{2} - \hat{S}_{2} \big|
=
2 
\cdot 
\mer,
\\
\sum_{\ub} \big(S_{\ub} - \hat{S}_{\ub} \big)^2
&=&
\big(S_{1} - \hat{S}_{1} \big)^2
+
\big(S_{2} - \hat{S}_{2} \big)^2
=
2 
\cdot 
\mer^2.
\end{eqnarray*}

One can state the result (\ref{eq: example_bound_max_3}) more formally in order to characterize the tightness of the bound presented in Theorem~\ref{th: si_error_general}.
\begin{theorem}
\label{th: si_error_general_lower_max}
For any $0 \leq \eps \leq 1$ there exists $f, \hat{f} \in L^2(\XX, \mu)$ such that $\V_\mu[f]~>~0$, $\V_\mu[\hat{f}] > 0$ and for corresponding Sobol' (total-effect) indices 
\begin{equation}
\label{eq: si_error_general_lower_max}
\max_{\ub  \subseteq \{1, \ldots, d\}}\Big|S_\ub - \hat{S}_\ub \Big|
= 
\eps 
\cdot 
\mer.
\end{equation}

\end{theorem}

Consequently, the error bound~(\ref{eq: si_error_general}) can be overestimated, but it is achievable.  Thus, the problem of accuracy of Sobol' indices estimates is reduced to the assessment of approximation quality.

%
%

\subsection{Risk bounds}

In this section, we study how Sobol' indices errors behave when the underlying approximation is constructed using the random design of experiments. At first, we generalize for this case the results of Section~\ref{sec: error_bounds}. Then, we consider the errors of Sobol' indices for two presented methods of approximation construction, depending on the size of the training sample and the complexity of  the approximating model.


\subsubsection{General results}


All risk bounds in this paper are obtained under the condition of random design:
\begin{Condition}
\label{cond: random_design}
Suppose the design $\des = \{\xb_i\}_{i=1}^{n}$ is randomly i.i.d. sampled from the probability measure $\mu$.
\end{Condition}

Consider an arbitrary approximation $\hat{f}$ of function $f$ constructed based on the training sample  $\mathcal{S}$   of size $n$. We assume that there is some fixed deterministic learning procedure that constructs the approximation from the given training sample (\ref{eq:training_sample}): 
\begin{equation}
\label{eq: appr_deterministic_training}
\mathcal{L}\colon \big(\des, \; Y=f(\des) + \etab \big) \to \hat{f},
\end{equation}
where fixed training sample leads to uniquely defined approximation $\hat{f}$.

The following theorem establishes a connection between the risk of estimated Sobol' indices $\E (S_\ub - \hat{S}_\ub)^2$ and $\E |S_\ub - \hat{S}_\ub |$, and the quadratic risk of the approximation $\E \|f - \hat{f}\|_{\mu}^2$  {\it if averaging occurs over random design points and noise realizations}. 

In what follows we will assume that $\V_{\mu}[f] > 0$. If we could guarantee that the approximation does not degenerate into constant, then the generalization of Theorem~\ref{th: si_error_general} to this setting would be straightforward. However, this is not the case, and we need to take into account the possibility of constant approximation, which leads to undefined Sobol' indices (see Definition~\ref{def: si}). For this reason, we additionally define\footnote{The choice of specific values is not  important in further conclusions.} Sobol' indices $\hat{S}_\ub = 2^{-d}$ in the case of $\V_{\mu} [\hat f] = 0$.

\begin{theorem}
\label{th: si_mae_risk_general}
Let $\hat{f}$ be an arbitrary  approximation of $f$  that  is constructed according to (\ref{eq: appr_deterministic_training}). Suppose that  under Condition~\ref{cond: random_design} of random design we have $\E \|f - \hat{f}\|_{\mu}^2 < \infty$. Then for corresponding Sobol' and total-effect indices of $f$ and $\hat{f}$ it holds for $\ub \subseteq \{1, \ldots, d\}$
\begin{eqnarray} 
\label{eq: mse_risk_1_exp_simple}
\max_{\ub} \, \E \big(S_\ub - \hat{S}_\ub \big)^2
&\leq&  
\ris^2,
\\
\label{eq: si_mae_general_each_index}
\E \big|S_\ub - \hat{S}_\ub \big|
&\leq&  
\ris  \left(\ris + 2\sqrt{S_{\ub}} \right),
\\
\label{eq: general_square_risk_normalized}
\text{where} \;\; 
\ris^2 
&\eq& 
\frac{\E{\|f - \hat{f}\|^2_{\mu} }}{\V_{\mu}[f]}.
\end{eqnarray}
\end{theorem}

\begin{corollary} 
\label{th: si_sum_mse_risk_general}
Under the assumptions of Theorem~\ref{th: si_mae_risk_general}, it holds for corresponding Sobol' indices\footnote{Not valid for total-effects.} for $\ub \subseteq \{1, \ldots, d\}$
\begin{eqnarray}
\label{eq: si_error_general_sum_squared_mse}
\E \Big[\sum_{\ub} \big(S_\ub - \hat{S}_\ub \big)^2 \Big]
\leq
2 \cdot \ris^2.
\end{eqnarray}
\end{corollary}

\begin{remark}
Note that Theorem~\ref{th: si_mae_risk_general} does not impose any  assumptions on the structure of approximating model, except (\ref{eq: appr_deterministic_training}). Particularly, the approximation does not need to have the tensor product form~(\ref{eq: hat_f_expansion}).
\end{remark}

\begin{remark}
Using Markov's inequality and (\ref{eq: mse_risk_1_exp_simple}, \ref{eq: si_mae_general_each_index}), one can obtain probabilistic bounds for (total) Sobol' indices estimate. Given $\gamma \in (0,1]$,
\begin{equation*}
P\Big\{ |S_\ub - \hat{S}_\ub | \geq \gamma \Big\}
\leq  
 \ris^2 / \gamma^2.
\end{equation*}
In addition, assume that Sobol' or total-effect index $S_\ub = 0$ for some $\ub$, then
\begin{equation*}
P\Big\{ \hat{S}_\ub  \geq \gamma \Big\}
\leq  
 \ris^2 / \gamma.
\end{equation*}

\end{remark}

\begin{remark}
Under the assumptions of Theorem~\ref{th: si_mae_risk_general}, if additionally $\E \|f - \hat{f}\|_{\mu}^4~<~\infty$, one can prove an estimate for $\E (S_\ub - \hat{S}_\ub)^2$ that is potentially more accurate for small indices:
\begin{equation*}
\E \big(S_\ub - \hat{S}_\ub \big)^2
\leq
2 \cdot \frac{\E{\|f - \hat{f}\|^4_{\mu} }}{\V^2_{\mu}[f]} 
+ 
8 S_\ub \cdot \frac{\E{\|f - \hat{f}\|^2_{\mu} }}{\V_{\mu}[f]}.
\end{equation*}
The proof repeats the proof of Theorem~\ref{th: si_mae_risk_general}, using $(\mer^2 + 2\sqrt{S_{\ub}} \cdot \mer)^2 \leq 2 \mer^4 + 8 S_{\ub} \cdot \mer^2$.
\end{remark}

\subsubsection{Projection method}

Theorem~\ref{th: si_mae_risk_general} considers an approximation of general type. Now, we will study the behavior of the particular method to construct the approximation, starting with the projection (quasi-regression) method defined as (\ref{eq: Projection_coeff}).

\begin{Condition}
\label{cond: f_bounded_XX}
We additionally require $f$ to be bounded on $\XX\colon$
\begin{equation}
\label{eq: f_bounded_L}
\big|f(\xb)\big| \leq L \;\;\; \text{for} \;\;\xb \in \XX.
\end{equation}
\end{Condition}

In this setting, it is straightforward to estimate a bound for the risk of the approximation (\ref{eq: general_square_risk_normalized}) and to obtain a corollary of Theorem~\ref{th: si_mae_risk_general} for the case of the projection method.
\begin{theorem}
\label{th: si_mse_int}
Under Condition~\ref{cond: random_design} of random design and boundedness Condition~\ref{cond: f_bounded_XX}, for corresponding Sobol' and total-effect indices of $f$ and $\hat{f}^P$ it holds for $\ub \subseteq \{1, \ldots, d\}$
\begin{eqnarray}
\label{eq: si_mse_int}
\max_{\ub} \, \E \big(S_\ub - \hat{S}_\ub^{P} \big)^2
&\leq&  
\ris_{p}^2,
\\
\label{eq: si_mae_int_each_index}
\E \big|S_\ub - \hat{S}_\ub^{P} \big|
&\leq&  
\ris_{p}  \left(\ris_{p} + 2 \sqrt{S_{\ub}} \right),
\\
\nonumber
\text{where} \;
\ris_{p}^2 
\eq
  \frac{1}{\V_{\mu}[f]}  &\cdot&
  \|e_N\|_{\mu}^2 
  + 
  \frac{L^2 + \sigma^2}{\V_{\mu}[f]}  \cdot \frac{N}{n}.
\end{eqnarray}

\end{theorem}


An obvious corollary of the previous theorem is the sufficient condition for the convergence of Sobol' indices in the mean square sense. 
\begin{corollary} 
Under the assumptions of Theorem~\ref{th: si_mse_int}, suppose additionally $\lim_{N \to \infty} \|e_{N}\|_{\mu} = 0$. Let $N = N(n)$,
\[
\frac{N}{n} \to 0 \;\; \text{and}  \;\; N \to \infty \;\; \text{as}  \;\; n \to \infty,
\] 
then
\[
\E \big(S_\ub - \hat{S}_\ub^{P} \big)^2 \ton 0.
\]
\end{corollary}

\subsubsection{Ordinary least squares}

The approximation of the least squares type~(\ref{eq: ls_empirical_minimization}) is especially sensitive to the design quality. Although according to (\ref{eq: ONB}),
\[
\E \left[\frac{\Phi^T \Phi}{n}\right] = I_N
\] 
($I_N$ is the identity matrix of size $N \times N$), there is some probability of ill-conditioned {\it normalized information matrix} $\Phi^T \Phi / n$ that may lead to unbounded approximation with 
\[
\hat{\cb} = \left(\frac{\Phi^T \Phi}{n} \right)^{-1} \frac{\Phi^T Y}{n}.
\]

As a result, typical risk bounds for the least squares method in the random design case are obtained for the approximation with additionally truncated   absolute values \cite{Kohler2002}. One consequence is that we cannot obtain risk bounds for Sobol' indices directly from Theorem~\ref{th: si_mae_risk_general}.






{\bf Least squares stability.} Let us introduce the numerical characteristic often used in the random design setting that helps us to control the design quality. 
\begin{definition}
For the orthonormal set $\{\Psi_{\al}, \; \al \in \rr_N\}$ that satisfies (\ref{eq: ONB}) and for some fixed nested sets of $d$-dimensional  multi-indices (\ref{eq: fixed_multiind_order}) define
\begin{equation}
\label{eq: K_N}
K_N \eq \sup_{\xb \in \XX} \Big[\sum_{\al \in\rr_N} \Psi_{\al}^2(\xb)\Big].
\end{equation}
\end{definition}


$K_N$ is connected with {\it the inverse Christoffel function} \cite{MR862231}, defined as
\[
k_N (\xb) \eq \sum_{\al \in\rr_N} \Psi_{\al}^2(\xb),
\] 
which is the diagonal of the orthogonal projection kernel onto $\pnu_N$, and
\[
K_N \eq \big\|k_N \big\|_{L^{\infty}}
=
\max_{u \in \pnu_N} \frac{\|u\|^2_{L^{\infty}}}{\|u\|^2_{\mu}}.
\]



$K_N$ depends on the measure $\mu$ and space $\pnu_N$. Although, the way one constructs the sequence (\ref{eq: fixed_multiind_order}) influences $K_N$, it is independent on the choice of the orthonormal basis in $\pnu_N$. It can be shown that $K_N \geq N$ and in general case $K_N$ may be arbitrary large~\cite{Cohen13}. In cases that we will consider in Section~\ref{sec: Interpretation}, $K_N \sim N$ or $\unsim N^2$. 

Denote {\it the spectral norm} of matrix $A \in \R^{m\times p}$ as
\[
|||A||| = \max_{\bz \in \R^p\colon\|\bz\| \neq 0}\frac{\|A\bz\|}{\|\bz\|}.
\]

Following \cite{Cohen13}, we state a lemma that allows to bound the spectral norm of the normalized information matrix\footnote{The proof can be  found in \cite{Cohen13}, see the corrected version.}.
\begin{lemma}
\label{le: spectral_norm_exponential_small}
Under Condition~\ref{cond: random_design} of random design, for $\delta \in (0, 1)$
\begin{equation}
P \Big\{ |||\Phi^T \Phi/n - I_N||| > \delta \Big\}  \leq 
2N \exp \left[-\frac{c_\delta \cdot n}{K_N} \right],
\end{equation}
where $c_\delta \eq (1+\delta)\ln(1+\delta)-\delta > 0$.
\end{lemma}

Lemma~\ref{le: spectral_norm_exponential_small} leads to the condition on $n$ and $N$ that excludes the possibility of ill-conditioned information matrix with high probability.
\begin{Condition}
\label{cond: ls_K_N_bound}
For some fixed $r > 0$ the relation of $N$ and $n$ satisfies
\begin{equation}
\label{eq: def_kappa_r_and_condition_KN_n}
K_N \leq \kappa_r \cdot \frac{n}{\ln n}, \;\; \text{where} \;\  \kappa_r = \frac{3 \cdot \ln(3/2) - 1}{2+2r}.
\end{equation}
\end{Condition}

Under the Condition~\ref{cond: ls_K_N_bound}, we have based on Lemma~\ref{le: spectral_norm_exponential_small} 
\begin{equation}
\label{eq: condition_K_N_and_n_05}
P \Big\{ |||\Phi^T \Phi/n - I_N||| > 1/2 \Big\}  \leq  2n^{-r}.
\end{equation}



{\bf Indices stability.} As the least squares estimate (\ref{eq: ls_coeff}) can be degenerate if $\det (\Phi^T \Phi) =  0$, we additionally define $\hat{c}_{\al} = 0$ in this case.

\begin{theorem}
\label{th: si_mse_ls_noise} 
Under Condition~\ref{cond: random_design} of random design and the stability Condition~\ref{cond: ls_K_N_bound}, for corresponding Sobol' and total-effect indices of $f$ and $\hat{f}^{LS}$ it holds for $\ub~\subseteq~\{1, \ldots, d\}$
\begin{eqnarray}
\label{eq: si_mse_ls_noise}
\max_{\ub} \, \E \big(S_\ub - \hat{S}_\ub^{LS} \big)^2
&\leq&  
\ris_{LS}^2 + 2n^{-r},
\\
\label{eq: si_mae_ls_each_index}
\E \big|S_\ub - \hat{S}_\ub^{LS} \big|
&\leq&
\ris_{LS}  \left(\ris_{LS} + 2 \sqrt{S_{\ub}} \right) + 2n^{-r},
\\
\nonumber
\text{where} \;\;
\ris_{LS}^2 
&\eq&
  \frac{1.2}{\V_{\mu}[f]}  \cdot
  \|e_N\|_{\mu}^2 
  + 
  \frac{4\sigma^2}{\V_{\mu}[f]}  \cdot \frac{N}{n}.
\end{eqnarray}

\end{theorem}

Note that, as Sobol' indices are bounded, we do not need to  truncate large values of the approximation additionally. Besides, the previous result does not require a maximum absolute value of the function.

Consider a noiseless case.
\begin{corollary}
\label{col: si_mse_ls}
Under the assumption of Theorem~\ref{th: si_mse_ls_noise} and provided noiseless observations i.e. $\sigma^2 = 0$, it holds
\begin{eqnarray}
\label{eq: si_mse_ls}
\E\big(S_\ub - \hat{S}_\ub^{LS} \big)^2
&\leq&  
   \frac{1.2}{\V_{\mu}[f]}  
  \|e_N\|_{\mu}^2 
  + 
  2n^{-r},
\\
\label{eq: si_mae_ls}
\E\big|S_\ub - \hat{S}_\ub^{LS} \big|
&\leq&  
   \frac{1.2}{\V_{\mu}[f]}  
  \|e_N\|_{\mu}^2 
  +
   \frac{2.2 \cdot S_{\ub}^{1/2}}{\V^{1/2}_{\mu}[f]}  
  \|e_N\|_{\mu} 
  + 
  2n^{-r}.
\end{eqnarray}


\end{corollary}

%


We have the following sufficient conditions for indices convergence in the mean square sense.
\begin{corollary} 
\label{th: converg_ls_asymp}
Under the assumptions of Theorem~\ref{th: si_mse_ls_noise} (except Condition~\ref{cond: ls_K_N_bound}), suppose additionally $\lim_{N \to \infty} \|e_{N}\|_{\mu} = 0$. Let $N = N(n)$,
\[
\frac{K_N \cdot \ln N}{n} \to 0 \;\; \text{and}  \;\; N \to \infty \;\; \text{as}  \;\; n \to \infty,
\] 
then
\[
\E\big(S_\ub - \hat{S}_\ub^{LS} \big)^2 \ton 0.
\]
\end{corollary}


\section{Asymptotic behavior}
\label{sec: Interpretation}

In this section, we discuss specific cases of approximation construction using various types of regressors, which are orthonormal w.r.t. different probability measures. For these particular approximations, we consider the interpretation of Theorem~\ref{th: si_mse_int} and \ref{th: si_mse_ls_noise} from an asymptotic point of view. The results of this section are summarized in Table~\ref{tab: SI_risk_all}.

We suppose that the analyzed function has given {\it smoothness} defined as below. 
\begin{definition}
\label{def: holder_condition}
A function $f$ on $\XX$ is said to satisfy a Holder condition with exponent $\beta \in (0, 1]$ if there is $\gamma > 0$  such that $|f(\xb) - f(\xb_0)| < \gamma \|\xb - \xb_0\|^{\beta}$ for $\xb, \xb_0 \in \XX$.
\end{definition}

\begin{definition}
\label{def: p-smoothness}
Let $m$ be a nonnegative integer, $\beta \in (0, 1]$. Set $p = m + \beta$. A function $f$ on $\XX$ is said to be $p$-smooth if it is $m$ times continuously differentiable on $\XX$ and 
\[
D^{\al}f = \frac{\partial^{|\al|}f}{\partial x_1^{\alpha_1}, \ldots, \partial x_d^{\alpha_d}}
\]
satisfies a Holder condition with exponent $\beta$ for all $\alpha$ with $|\al| = \sum_{i=1}^d \alpha_i = m$.
\end{definition}

Thus, the function under analysis $f(\xb)$ is assumed to be $p$-smooth.

\subsection{Ordinary least squares in the noiseless case}
\label{sec: PCE_asymp_int}

Following Corollary~\ref{col: si_mse_ls}, we start with the approximations based on OLS in the~case of noiseless observations ($\sigma^2 = 0$, see the observation model in Section~\ref{sec: Appr_constr}).

\subsubsection{Legendre polynomials}
\label{sec: legendre_asymp_int} 

Consider a uniform input distribution on $\XX = [-1, 1]^d$. We exploit Legendre polynomials and ordinary least squares to approximate $f(\xb)$ and then calculate  Sobol' indices based on this expansion. The regressors are constructed as a tensor product of normalized univariate Legendre polynomials: 
\begin{equation}
\label{eq: legendre_multivariate}
\Big\{\Psi_{\al}(\xb) = \prod_{i=1}^d \psi_{\alpha_i}(x_i), \; \al \in \rr_N \Big\},
\end{equation}
where  $\psi_{\alpha}(x) = \tilde{\psi}_{\alpha}(x) /\|\tilde{\psi}_{\alpha}\|_{\mu}$ and $\tilde{\psi}_{\alpha}(x)$ are univariate Legendre polynomials:
\begin{eqnarray*}
\tilde\psi_{0} &=& 1, \;\; 
\tilde\psi_{1} = x, \;\;
\tilde\psi_{2} = \frac{1}{2}(3x^2-1),  \;\;\tilde\psi_{3} = \frac{1}{2}(5x^3 - 3x), \;\; \ldots,
\\
\tilde\psi_{\alpha} 
&=& 
\frac{1}{2^{\alpha}\alpha!}\frac{d^{\alpha}}{dx^{\alpha}} (x^2-1)^{\alpha} 
\;\;
\text{with}
\;\;
\|\tilde\psi_{\alpha}\|_{\mu} = (2\alpha +1)^{-1/2}, \;\; \alpha \geq 0.
\end{eqnarray*} 

In all cases in this section to construct the  truncation set $\rr_N$, we use {\it maximum degree} truncation scheme that is suitable for our asymptotic analysis. For some $\alpha_{max} \in \NN_{+}$ called maximum degree we have
\begin{equation}
\label{eq: max_degree}
\rr_N = \{\al \in \NN^d\colon  \max_{i=1,\ldots,d} \{\alpha_i\} \leq \alpha_{max}\},
\end{equation}
with  $N = |\rr_N| =  (\alpha_{max}+1)^d$.

After selection of the probability measure and the regressors, we can estimate $K_N$ to ensure the stability Condition~\ref{cond: ls_K_N_bound}:
\begin{eqnarray*}
K_N  
&=& 
\sup_{\xb \in \XX} \Big[\sum_{\al \in\rr_N} \Psi_{\al}^2(\xb)\Big]
=
\sup_{\xb \in [-1,1]^d} \Big[\sum_{\al \in\rr_N} \prod_{i=1}^d \psi_{\alpha_i}^2(x_i)\Big]
\\
&=& 
\sum_{\al \in\rr_N} \prod_{i=1}^d \psi_{\alpha_i}^2(1)
= 
\prod_{i=1}^d  \sum_{\alpha = 0}^{\alpha_{max}} (2\alpha+1) = (\alpha_{max} +1)^{2d} = N^2,
\end{eqnarray*}
here we used $\|\psi_{\alpha}\|_{L^{\infty}([-1,1])}  \eq \sup_{[-1,1]} |\psi_{\alpha}(x)| = \psi_{\alpha}(1) = (2\alpha+1)^{1/2}$. Calculation of $K_N$ for more general truncation schemes can be found in \cite{Chkifa15}.

Thus, the stability Condition~\ref{cond: ls_K_N_bound} takes the form:
\begin{equation}
\label{eq: N_opt_ls_leg_noise-free}
N \leq \kappa_r^{1/2} \cdot \left[\frac{n}{\ln n}\right]^{1/2}.
\end{equation}

It remains to estimate the error $\|e_N\|_{\mu}$.  One can do it using the inequality
\[
\|e_N\|_{\mu}=\inf_{g \in \pnu_N}\|f-g\|_{\mu}  \leq \inf_{g \in \pnu_N}\|f-g\|_{L^\infty(\XX)}.
\] 
Denote $ d_N(f) \eq \inf_{g \in \pnu_N}\|f-g\|_{L^\infty(\XX)}$. Asymptotic bounds for $d_N(f)$ are obtained in classical statistical literature. The space $\pnu_N = span\{\Psi_{\al}, \; \al \in \rr_N\}$ of multivariate Legendre polynomials with fixed maximal degree coincides with the space of multivariate algebraic polynomials with the same maximal degree\footnote{The basis of this space can be constructed as a tensor product of univariate algebraic polynomials $1$, $x$, $x^2$, \ldots, $x^{\alpha_{max}}$.}, therefore  this two spaces have the same $d_N(f)$.  According to \cite{Huang98}, for multivariate algebraic polynomials of maximal degree $\alpha_{max}$ and $p$-smooth function $f$ we have up to constant factor\footnote{Given two sequences of positive numbers $\{a_n\}$ and $\{b_n\}$, the notation $a_n \lesssim b_n $ means that $a_n / b_n$ is bounded.} $
d_N(f) \lesssim \alpha_{max}^{-p}$ and  $d_N(f) \lesssim N^{-p/d}$.


Combining all parts, based on Corollary~\ref{col: si_mse_ls}, we obtain an inequality
\[
\E\big(S_\ub - \hat{S}_\ub^{LS} \big)^2 \lesssim 
\|e_N\|_{\mu}^2 
  + 
  n^{-r}
\lesssim 
N^{-2p/d} + n^{-r}
\lesssim 
\left(\frac{n}{\ln n}\right)^{-p/d} + n^{-r}.
\]
Setting $r = p/d$, we have for Legendre polynomials
\begin{equation}
\label{eq: si_ols_asymptotic_legendre}
\E\big(S_\ub - \hat{S}_\ub^{LS} \big)^2 \lesssim 
\left(\frac{n}{\ln n}\right)^{-p/d}.
\end{equation}

For comparison see Table~\ref{tab: SI_risk_all}.

\subsubsection{Chebyshev polynomials}

Similarly to previous section, we continue our analysis for Chebyshev polynomials and the arcsine input distribution on $\XX = [-1, 1]^d$ having the density
\[
p(\xb) = \prod_{i=1}^d \Big(\pi\sqrt{1-x_i^2} \Big)^{-1}, \;\; \xb \in [-1,1]^d.
\]

The regressors are constructed as a tensor product of normalized univariate Chebyshev polynomials: 
\[
\psi_0 = 1, \;
\psi_1 = \sqrt{2} \cdot x, \; 
\psi_2 = \sqrt{2} \cdot  (2x^2-1), \;
\ldots, 
\psi_\alpha = \sqrt{2} \cdot  \cos(\alpha \arccos x),
\]
with $K_N$ estimate:
\begin{eqnarray*}
K_N  
&=& 
\sup_{\xb \in \XX} \Big[\sum_{\al \in\rr_N} \Psi_{\al}^2(\xb)\Big]
=
\sup_{\xb \in [-1,1]^d} \Big[\sum_{\al \in\rr_N} \prod_{i=1}^d \psi_{\alpha_i}^2(x_i)\Big]
\\
&=& 
\sum_{\al \in\rr_N} \prod_{i=1}^d \psi_{\alpha_i}^2(1)
\leq
\prod_{i=1}^d  \sum_{\alpha = 0}^{\alpha_{max}} 2 = 2^d \cdot (\alpha_{max} +1)^{d} = 2^d \cdot N,
\end{eqnarray*}
we used $\|\psi_{\alpha}\|_{L^{\infty}([-1,1])} = \psi_{\alpha}(1) \leq \sqrt{2}$. The stability Condition~\ref{cond: ls_K_N_bound} takes the form
\begin{equation}
\label{eq: si_ols_asymptotic_cheb}
N \leq \kappa_r \cdot 2^{-d} \cdot \frac{n}{\ln n} .
\end{equation}

The error $\|e_N\|_{\mu} \leq d_N(f)$ is estimated similarly to Section~\ref{sec: legendre_asymp_int}.  Thus,
\[
\E\big(S_\ub - \hat{S}_\ub^{LS} \big)^2 \lesssim 
\|e_N\|_{\mu}^2 
  + 
  n^{-r}
\lesssim 
N^{-2p/d} + n^{-r}
\lesssim 
\left(\frac{n}{\ln n}\right)^{-2p/d} + n^{-r}.
\]
Setting $r = 2p/d$, we obtain for Chebyshev polynomials
\begin{equation}
\label{eq: si_ols_asymptotic_chebyshev}
\E\big(S_\ub - \hat{S}_\ub^{LS} \big)^2 \lesssim 
\left(\frac{n}{\ln n}\right)^{-2p/d}.
\end{equation}

\subsubsection{Trigonometric polynomials}

We proceed with the analysis for a uniform input distribution defined on $\XX = [0, 1]^d$ and Trigonometric polynomials. For this case, we additionally require that $f(\xb)$ can be extended to a $1$-periodic function in each of its arguments (which is equivalent to additional boundary conditions).

The regressors are constructed as a tensor product of normalized univariate Trigonometric polynomials: 
\begin{equation}
\begin{aligned}
\label{eq: Trigonometric}
\psi_0 &= 1, \;
\psi_1 = \sqrt{2} \cdot  \sin 2\pi x, \; 
\psi_2 = \sqrt{2} \cdot  \cos 2\pi x, \;
\\
\psi_3 &= \sqrt{2} \cdot  \sin 4\pi x, \; 
\psi_4 = \sqrt{2} \cdot  \cos 4\pi x, \;
\ldots
\end{aligned}
\end{equation}

Assume $\alpha_{max}$ is even. We have the following $K_N$ estimate:
\begin{eqnarray*}
K_N 
&=& 
\sup_{[0,1]^d} \Big[\sum_{\al \in\rr_N} \Psi_{\al}^2(\xb)\Big] 
\\
&=& \sup_{[0,1]^d} \Big[ \big(1 +2\sin^2(2\pi x_1) + 2\cos^2(2\pi x_1) + \ldots +2\cos^2(\alpha_{max}\pi x_1) \big)
\ldots
\\
&\ldots& \big(1 +2\sin^2(2\pi x_d) + 2\cos^2(2\pi x_d) +\ldots +2\cos^2(\alpha_{max}\pi x_d) \big)\Big]   
\\
&=& 
\left(1 + 2 \frac{\alpha_{max}}{2}\right)^d = N.
\end{eqnarray*}

The stability Condition~\ref{cond: ls_K_N_bound} takes the form
\begin{equation}
\label{eq: si_ols_asymptotic_trig}
N \leq \kappa_r \cdot \frac{n}{\ln n}.
\end{equation}

The error $\|e_N\|_{\mu} \leq d_N(f)$. According to \cite{Huang98,Rafajlowicz88}, for such $p$-smooth functions $d_N(f) \lesssim (\alpha_{max}/2)^{-p} \lesssim  N^{-p/d}$.  Thus,
\[
\E\big(S_\ub - \hat{S}_\ub^{LS} \big)^2 \lesssim 
\|e_N\|_{\mu}^2 
  + 
  n^{-r}
\lesssim 
N^{-2p/d} + n^{-r}
\lesssim 
\left(\frac{n}{\ln n}\right)^{-2p/d} + n^{-r}.
\]
Setting $r = 2p/d$, we obtain for Trigonometric polynomials
\begin{equation}
\label{eq: si_ols_asymptotic_trigonometric}
\E\big(S_\ub - \hat{S}_\ub^{LS} \big)^2 \lesssim 
\left(\frac{n}{\ln n}\right)^{-2p/d}.
\end{equation}

\subsection{Noisy case}

In the case of observations with noise, we have similar  results for Legendre, Chebyshev, and Trigonometric polynomials. Indeed, for the projection method of  coefficients estimation, following Theorem~\ref{th: si_mse_int} obtain
\[
\E\big(S_\ub - \hat{S}_\ub^{P} \big)^2
\lesssim 
  \|e_N\|_{\mu}^2 
  + 
  \frac{N}{n}
\lesssim
  N^{-2p/d}
  + 
  \frac{N}{n}.
\]
Balancing the summands at the right side, find asymptotically optimal $N(n)$ that minimizes the quadratic risk:
\begin{equation}
\label{eq: N_opt_proj_noise}
N = n^{\frac{d}{2p+d}},
\end{equation}
and the corresponding quadratic risk
\begin{equation}
\label{eq: integ_noise_assymp}
\E\big(S_\ub - \hat{S}_\ub^{P} \big)^2
\lesssim 
n^{-\frac{2p}{2p+d}}.
\end{equation}
The presented asymptotics is not improved even if $\sigma^2=0$. This rate corresponds to the Stone's bound \cite{Stone82}.

Back to ordinary least squares. Based on Theorem~\ref{th: si_mse_ls_noise} for Legendre, Chebyshev, and Trigonometric polynomials, it holds
\[
\E\big(S_\ub - \hat{S}_\ub^{LS} \big)^2
\lesssim 
   N^{-2p/d}
  + 
  \frac{N}{n}
  + n^{-r},
  \;\;
  \text{with}
  \;\;
  K_N \lesssim \frac{n}{\ln n}.
\]
Balancing the summands, obtain asymptotically optimal parameters that minimizes the risk (without the stability restriction on $K_N$)
\begin{equation}
\label{eq: N_opt_ls_noise}
N = n^{\frac{d}{2p+d}},
\;\;
r = 2p / (2p+d).
\end{equation} 

The risk for Chebyshev and Trigonometric polynomials:
\begin{equation}
\label{eq: ls_noise_assymp}
\E\big(S_\ub - \hat{S}_\ub^{LS} \big)^2
\lesssim 
n^{-\frac{2p}{2p+d}},
\end{equation}
which corresponds to the Stone's bound. Finally for Legendre polynomials, due to the restriction $K_N \lesssim n / \ln n$, we have 
\begin{equation*}
  \E\big(S_\ub - \hat{S}_\ub^{LS} \big)^2
\lesssim 
 \begin{cases}
    n^{-\frac{2p}{2p+d}}, & \text{if $p/d > 1/2$};\\
    \left(\frac{n}{\ln n}\right)^{-p/d}, & \text{if $p/d \leq 1/2$}.
  \end{cases}
\end{equation*}

Table~\ref{tab: SI_risk_all} summarizes asymptotic bounds for the risk of Sobol' indices obtained in a different setting.

\begin{remark}
Although the obtained noisy rate (\ref{eq: ls_noise_assymp}) for the least squares method coincides with the Stone's bound, the corresponding noiseless rates (\ref{eq: si_ols_asymptotic_chebyshev}, \ref{eq: si_ols_asymptotic_trigonometric}) outperform this bound. There is no contradiction here, because Stone's assumptions \cite{Stone82} imply $\V(y|\xb) > 0$. See \cite{BAUER201793,KOHLER2014197,KOHLER20131871} for more information on noiseless rates.
\end{remark}


\section{Experiments}
\label{sec: Experiments}

\subsection{Test functions}

The following two functions are commonly used for benchmarking in global sensitivity analysis.

{\bf Sobol' g-function.} Let $\xb$ be uniformly distributed in the hypercube $[0, 1]^d$,
\[
f(\xb) = \prod_{i=1}^d \frac{|4x_i-2|+c_i}{1+c_i}, \; \text{with} \; c_i \geq 0, \; i = 1, \ldots, d.
\]
In our case $d=2$, and we use two sets of parameters in different experiments:  $\{c_1 = 0, \; c_2=4\}$ and $\{c_1 = 0, \; c_2=0.5\}$. Analytical expressions for the corresponding Sobol' indices are available in \cite{Sobol93}:
\[
S_{\ub} = \frac{1}{D} \prod_{j \in \ub}\frac{1}{3}(c_j +1 )^{-2},  \;\; \ub \subseteq \{1, \ldots, d\},
\]
where $D = \prod_{i=1}^d \big[1 + \frac{1}{3} (c_i +1)^{-2} \big] - 1$.

{\bf Ishigami function.} Let $\xb$ be uniformly distributed in the hypercube $[-\pi, \pi]^3$,
\[
f(\xb) = \sin x_1 + a \cdot  \sin^2 x_2 + b x_3^4 \cdot \sin x_1, \;\; a=7, \; b=0.1.
\]
Theoretical values for its Sobol' indices can be found in \cite{Saltelli2008}:
\begin{eqnarray*}
S_1 &=&  \frac{b \pi^4}{5v} + \frac{b^2 \pi^8}{50v} +\frac{1}{2v}, \;\;
S_2 = \frac{a^2}{8v}, \;\;
S_3 = 0,\;\;
\\
S_{1 3} &=& \frac{8 b^2 \pi^8}{225v},\;\;
S_{1 2} = S_{2 3} = S_{1 2 3} = 0,
\end{eqnarray*}
where $v = a^2/8 + b \pi^4/5 + b^2 \pi^8 / 18 + 1/2$.

\subsection{Error bounds}
\label{sec: experiments_error_bounds}

Experiments in this section illustrate the results provided in \ref{sec: error_bounds}, including justification of theoretical error bounds and the demonstration of the method proposed in \ref{sec: a_new_method}. Observations are supposed to be noiseless ($\sigma^2 = 0$).

{\bf Approximating function.} All results of  Section~\ref{sec: error_bounds} are valid for arbitrary types of approximations. As an example, we use an approximation of the test functions based on multivariate Legendre polynomials (\ref{eq: legendre_multivariate}), the least squares method, and {\it hyperbolic truncation scheme} \cite{BlatmanSudret10b} corresponding to
\[
\rr_N = \big\{\al \in \NN^d\colon \|\al\|_q\ \leq t \big\}, \;
\]
where $ \|\al \|_q \triangleq \big( \sum_{i=1}^d \alpha_i^q \big)^{1/q}$ with fixed $q \in (0, 1]$ and $t \in \NN_+$ is a fixed maximal total degree of polynomials. The parameters of the truncation scheme are the following:
\begin{itemize}
\item {Sobol' g-function:} $q=0.5$, $t=20$, $N=91$ regressors;
\item { Ishigami function:} $q=0.75$, $t=20$, $N=815$ regressors;
\end{itemize}

\begin{figure}[t!]
    \centering
        \subfloat[Results for Sobol' g-function]{\includegraphics[width=\textwidth]{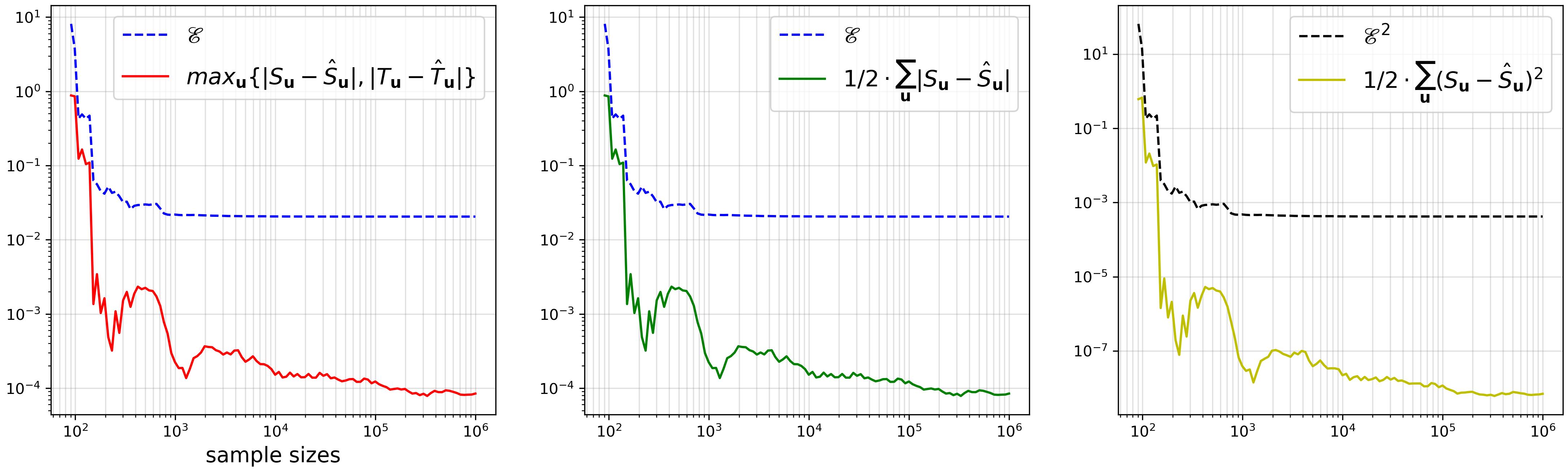}}
    \label{fig:subfigname_1}
        \qquad
        \subfloat[Results for Ishigami function\label{fig: theoretical_max_and_sum_bounds_Ishigami}]{\includegraphics[width=\textwidth]{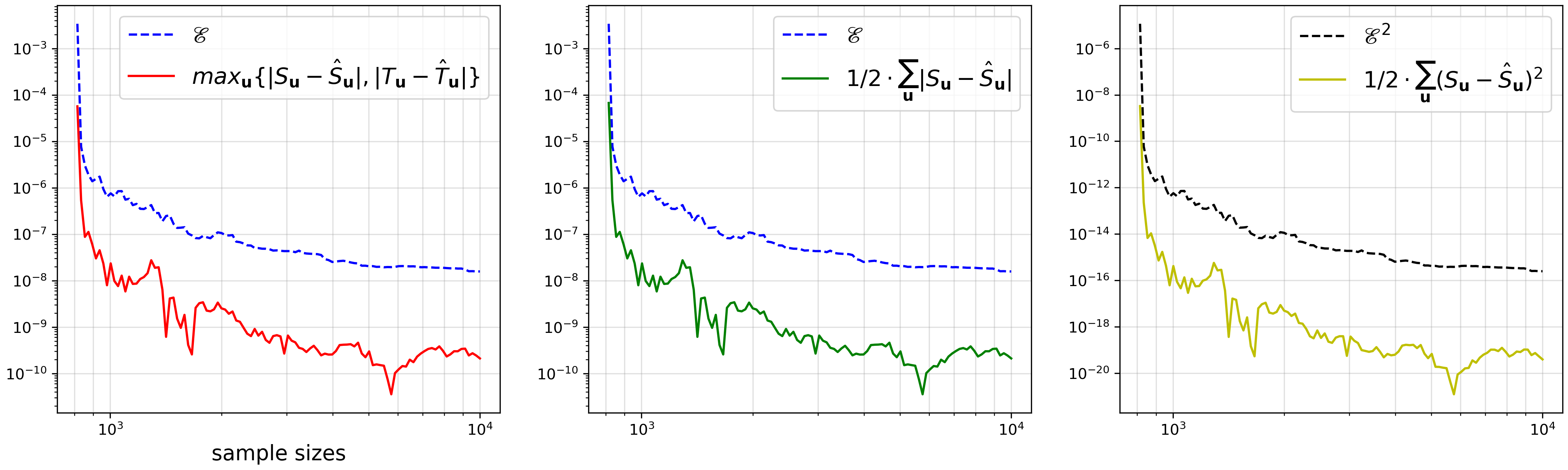}}
\caption{Different error metrics for Sobol' and total-effect indices and their theoretical error bounds (\ref{eq: si_error_general}, \ref{eq: si_error_general_sum_abs}, \ref{eq: si_error_general_sum_squared}). The approximation error $\mer$ is given by~(\ref{eq: def_relative_eps}). Fixed  number of regressors.}
\label{fig: theoretical_max_and_sum_bounds}
\end{figure}

\begin{figure}[t!]
    \centering
        \subfloat[Theoretical bounds: (\ref{eq: si_error_general}) and $\text{bound}(\mer, S_\ub) \;/ \; \text{bound}(\mer, T_\ub)$ given by (\ref{eq: si_error_general_each_S_approx})\label{fig: Sobol_theoretical_each}]{\includegraphics[width=\textwidth]{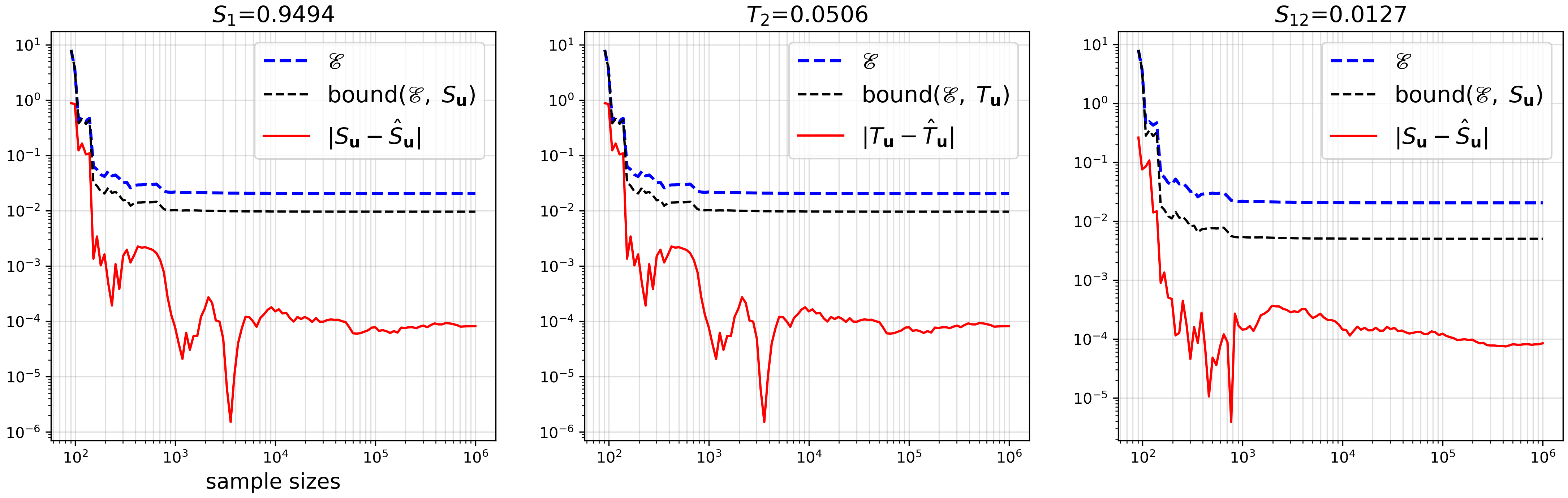}}
        \qquad
        \subfloat[Sample-based bounds: method proposed in \ref{sec: a_new_method} and  bootstrap (\ref{eq: bootstrap_appr_bound})\label{fig: Sobol_practical_each}]{\includegraphics[width=\textwidth]{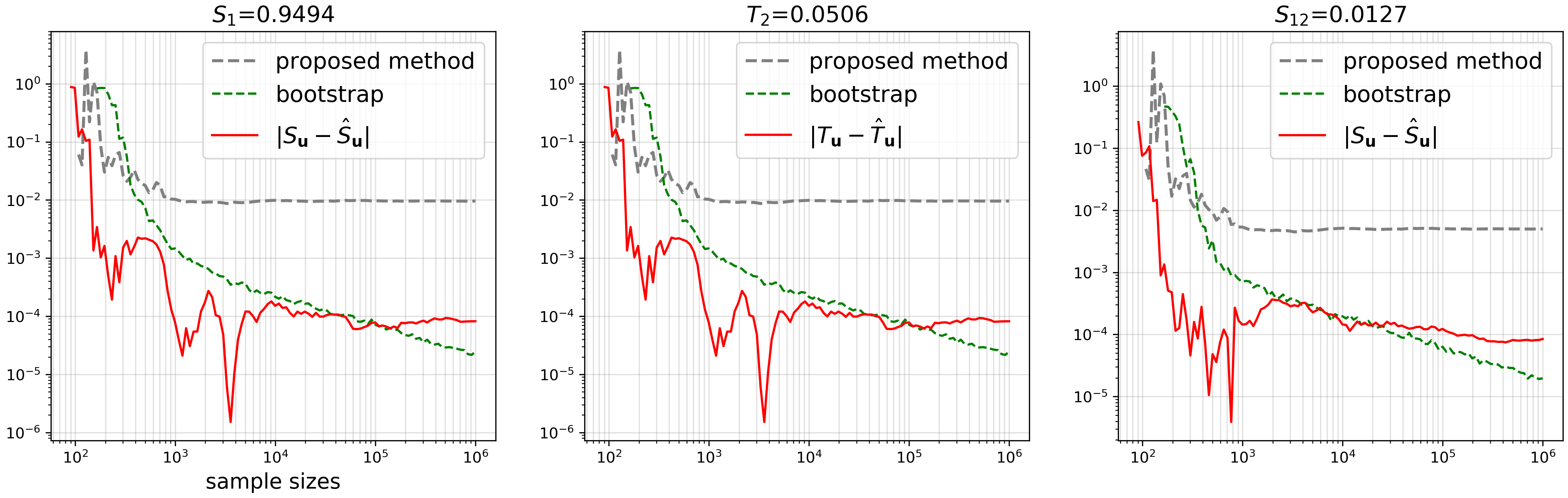}}
\caption{Sobol' g-function. Errors of selected Sobol' indices $|S_\ub - \hat{S}_\ub |$ and total-effects $|T_\ub - \hat{T}_\ub |$, and their theoretical and sample-based error bounds. True indices values of $S_1$, $T_{2}$ and $S_{12}$ are provided. Fixed $N=91$ regressors.}
\end{figure}

\begin{figure}[t!]
    \centering
        \subfloat[Theoretical bounds:  (\ref{eq: si_error_general}, \ref{eq: si_error_general_each_S_approx}, \ref{eq: si_error_general_zero-one}), where $\text{bound}(\mer, \cdot) =  \mer^2$ for $T_{12}$ and $S_{123}$ \label{fig: Ishigami_theoretical_each}]{\includegraphics[width=\textwidth]{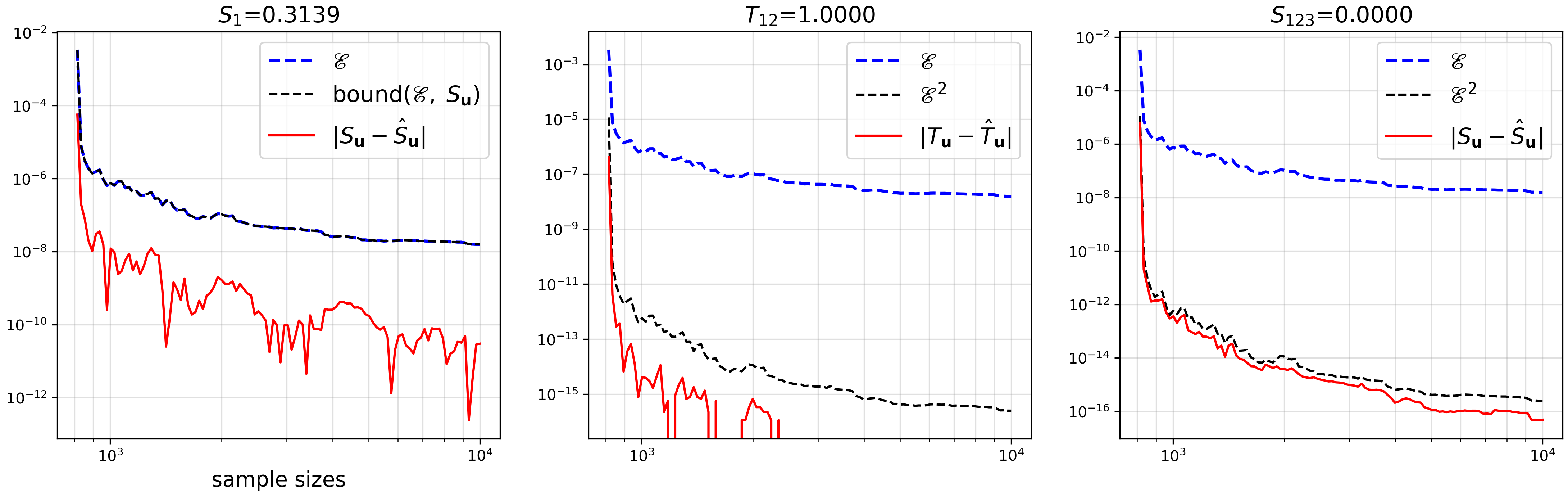}}
        \qquad
        \subfloat[Sample-based bounds: method proposed in \ref{sec: a_new_method} and  bootstrap (\ref{eq: bootstrap_appr_bound})\label{fig: Ishigami_practical_each}]{\includegraphics[width=\textwidth]{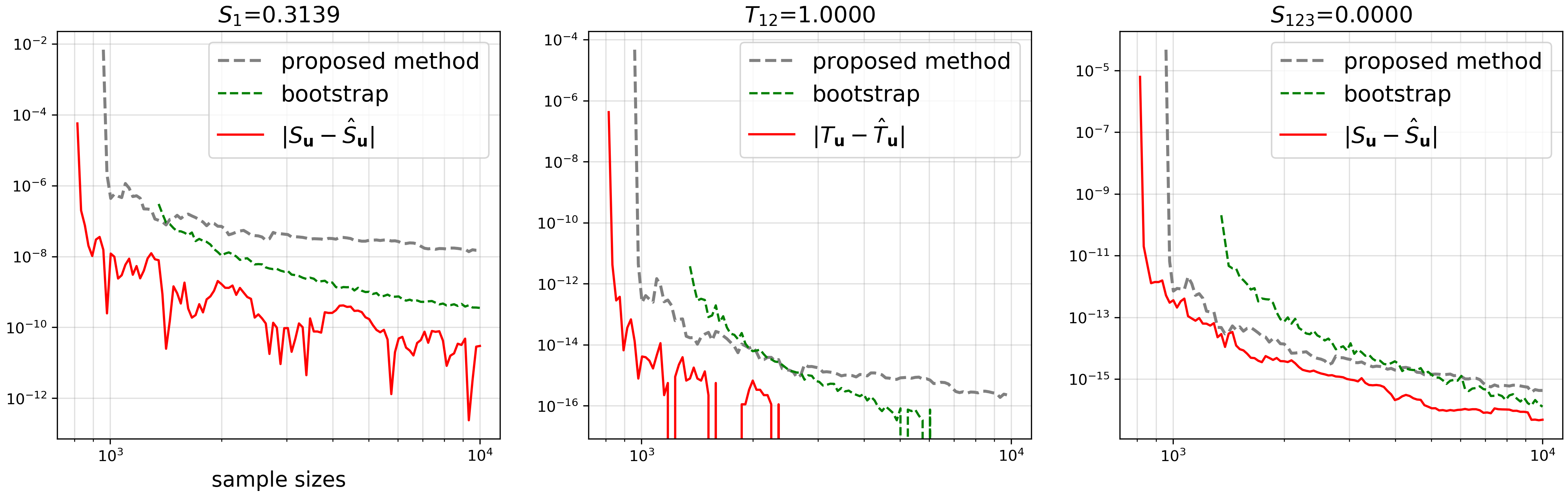}}
\caption{Ishigami function. Errors of selected Sobol' indices $|S_\ub - \hat{S}_\ub |$ and total-effects $|T_\ub - \hat{T}_\ub |$, and their theoretical and sample-based error bounds. True indices values of $S_1$, $T_{12}$ and $S_{123}$ are provided. Fixed $N=815$ regressors.}
\end{figure}

{\bf Theoretical error bounds} are shown in Figures~\ref{fig: theoretical_max_and_sum_bounds}, \ref{fig: Sobol_theoretical_each}, \ref{fig: Ishigami_theoretical_each}. This group of experiments illustrates Theorem~\ref{th: si_error_general}, Corollaries~\ref{th: si_error_general_resolved} and \ref{th: si_error_general_max}.  Particularly, bounds (\ref{eq: si_error_general}, \ref{eq: si_error_general_each_S_approx}, \ref{eq: si_error_general_zero-one}, \ref{eq: si_error_general_sum_abs}, \ref{eq: si_error_general_sum_squared}) are demonstrated. In these figures, ``$\text{bound}(\mer, S_\ub)$'' and ``$\text{bound}(\mer, T_\ub)$'' correspond to the error bound (\ref{eq: si_error_general_each_S_approx}), and may coincide with the relative error $\mer$ or its square $\mer^2$. 

The relative error of approximation $\mer$ is estimated based on (\ref{eq: MSE_noiseless_app}) with an independent test sample of size $10^6$. As the regressors number is fixed, the quality of the model and indices accuracy cease to improve, starting with some sample size.

As one can see, most of the presented theoretical bounds overestimate the corresponding errors of Sobol' indices, particularly for large sample sizes. Another observation is that estimates of small indices do tend to converge faster. Besides, the bound (\ref{eq: si_error_general_each_S_approx}) is more accurate for such indices and especially good for zero Sobol' and total-effect indices (in these cases for Ishigami function this bound coincides with $\mer^2$). The same holds for the indices close to $1$.

{\bf Sample-based error bounds} are shown in Figures~\ref{fig: Sobol_practical_each}, \ref{fig: Ishigami_practical_each}. These bounds are based on relatively small samples and can be used in practice.

%


{\it The proposed method} allows to estimate the error of Sobol' and total-effect indices obtained from the approximating function. It is described in \ref{sec: a_new_method} and is based on (\ref{eq: si_error_general_each_S_practice}). Following (\ref{eq: MSE_noiseless_app}), we use $15\%$ of the sample for holdout control of the approximation error.

For comparison, we examine sample-based error bounds based on the {\it bootstrap method} motivated by \cite{CIbootstrapSobol}. It is used $n_{s} = 100$ bootstrap subsamples  to get error bounds for estimated Sobol' indices from  the formula of the sample standard deviation
\begin{equation}
\label{eq: bootstrap_appr_bound}
err_\ub
=
3 \cdot \sqrt{\frac{1}{n_{s}} 
\sum_{j=1}^{n_s} 
\big(
\hat{S}^{(j)}_{\ub} - \bar{S}_{\ub}
\big)^2}, \;\; 
\text{with} \;\;
\bar{S}_{\ub} \eq \frac{1}{n_s}  \sum_{j=1}^{n_s} \hat{S}^{(j)}_{\ub},
\end{equation} 
where $\hat{S}^{(j)}_{\ub}$ is the Sobol' index of the subset $\ub \subseteq \{1, \ldots, d\}$  of input variables obtained from the approximation described above, which was constructed based on the $j$-th bootstrap subsample. The  error bounds for total-effects are calculated similarly. Note that due to random sampling with replacement and the lack of regularization in the least squares method, small sample sizes may lead to degenerate approximations and undefined bootstrap error bounds.

The experiments demonstrate that the bootstrap bounds may be overpessimistic for small sample sizes and overoptimistic for large samples. The latter can be explained by the fact that the method neglects the bias of the approximation, and its error estimate for indices can converge to zero even for inaccurate~approximations.

Compared to bootstrap, the proposed method gives more conservative error estimates for large samples and relatively accurate estimates for small ones. Another of its properties is possible instability for very small samples, which is a consequence of high uncertainty of the approximation error obtained with holdout control. As it should be expected, the corresponding sample-based bound converges to the theoretical bound (\ref{eq: si_error_general_each_S_approx}). Similar to the bound (\ref{eq: si_error_general_each_S_approx}), the proposed method is more accurate for small (and close to $1$) Sobol' and total-effect indices.

\subsection{Risk bounds}

\begin{figure}[t!]
    \centering
    
        \subfloat[Least squares method for various regressors numbers. Noiseless case, $\sigma=0$. \label{fig: Sobol_risk_ls}]{\includegraphics[width=\textwidth]{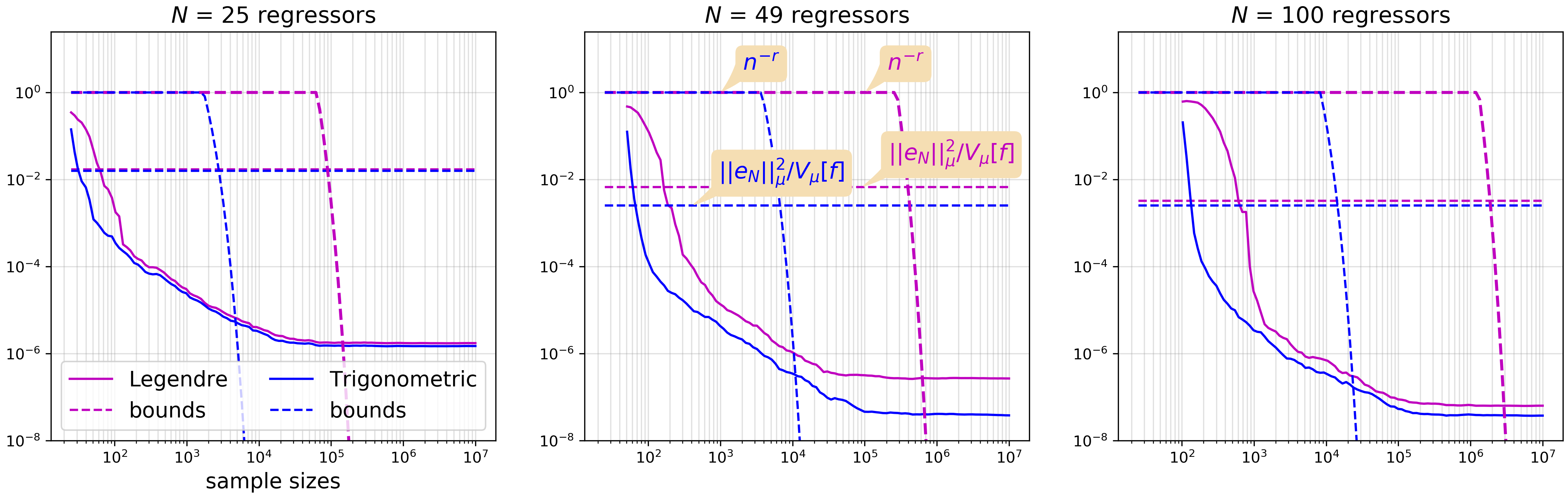}}
        \qquad
        \subfloat[Projection method for various regressors numbers. Noiseless case, $\sigma=0$. \label{fig: Sobol_risk_proj}]{\includegraphics[width=\textwidth]{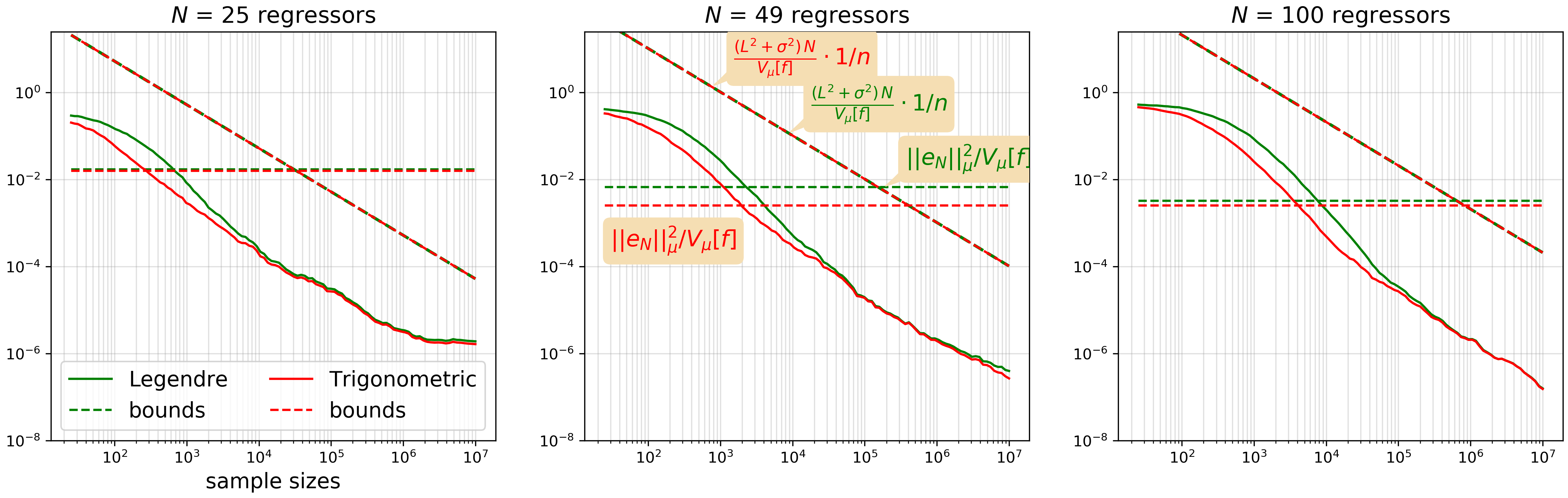}}
            \qquad
            \subfloat[Two methods for various noise levels $\sigma$. Legendre polynomials. $N=25$ regressors. \label{fig: Sobol_risk_noise}]{\includegraphics[width=\textwidth]{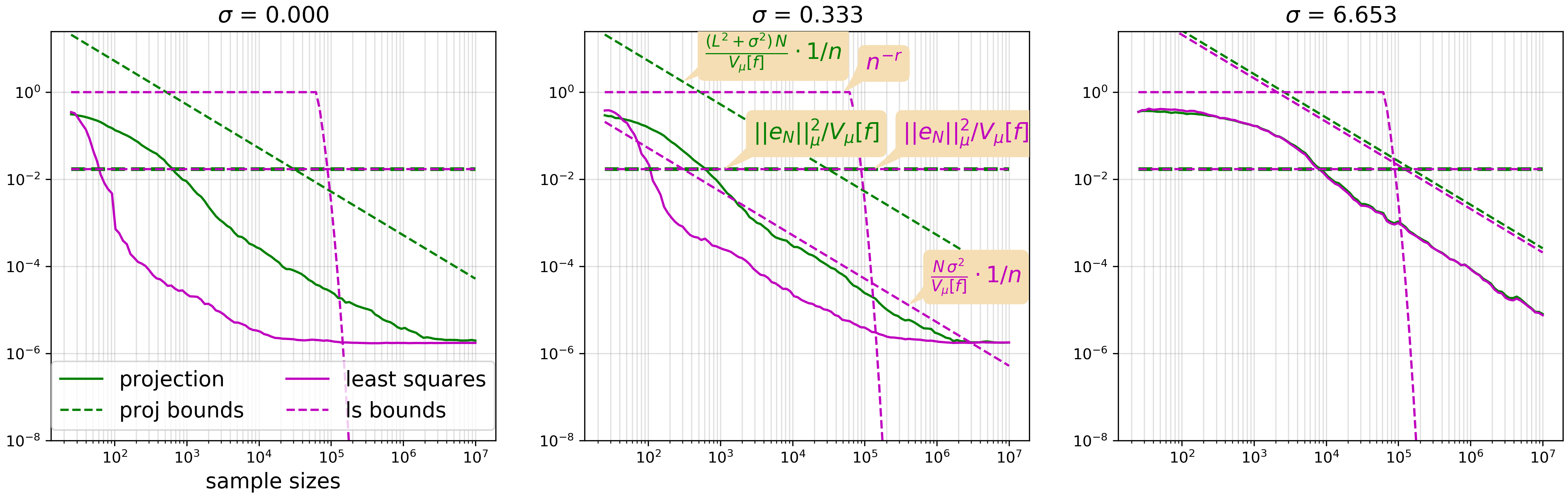}}
\caption{Sobol' g-function.  Quadratic risk of Sobol' and total-effect indices $\max_{\ub} \big\{\E (S_\ub - \hat{S}_\ub)^2, \; \E (T_\ub - \hat{T}_\ub )^2 \big\}$ and the components of bounds (\ref{eq: si_mse_int}, \ref{eq: si_mse_ls_noise}), depending on the regressors number and noise level.}
\end{figure}

\begin{figure}[t!]
    \centering
        \subfloat[Least squares method for various regressors numbers. Noiseless case, $\sigma=0$. \label{fig: Ishigami_risk_ls}]{\includegraphics[width=\textwidth]{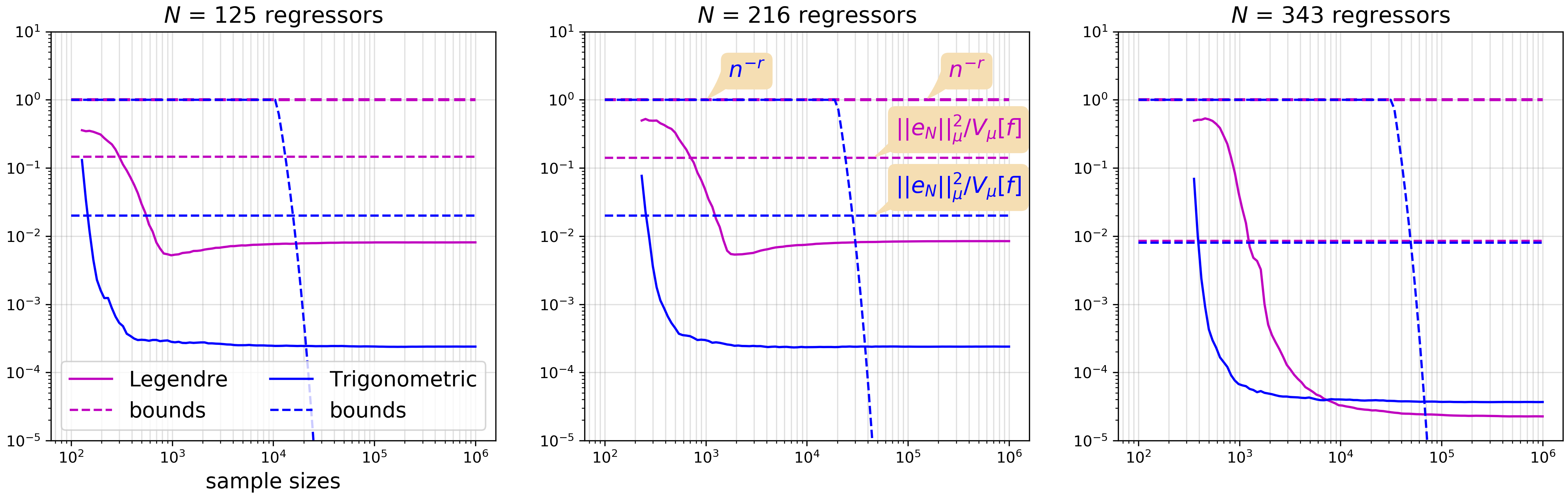}}
        \qquad
        \subfloat[Projection method for various regressors numbers. Noiseless case, $\sigma=0$. \label{fig: Ishigami_risk_proj}]{\includegraphics[width=\textwidth]{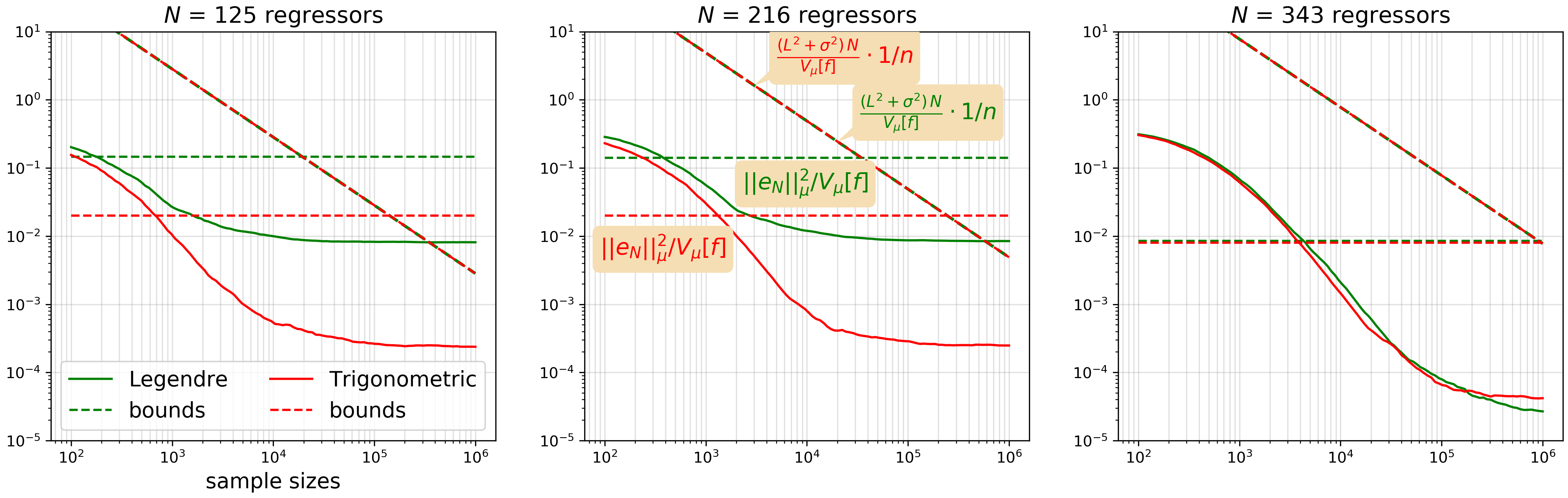}}
                \qquad
                \subfloat[Two methods for various noise levels $\sigma$. Legendre polynomials. $N=216$ regressors. \label{fig: Ishigami_risk_noise}]{\includegraphics[width=\textwidth]{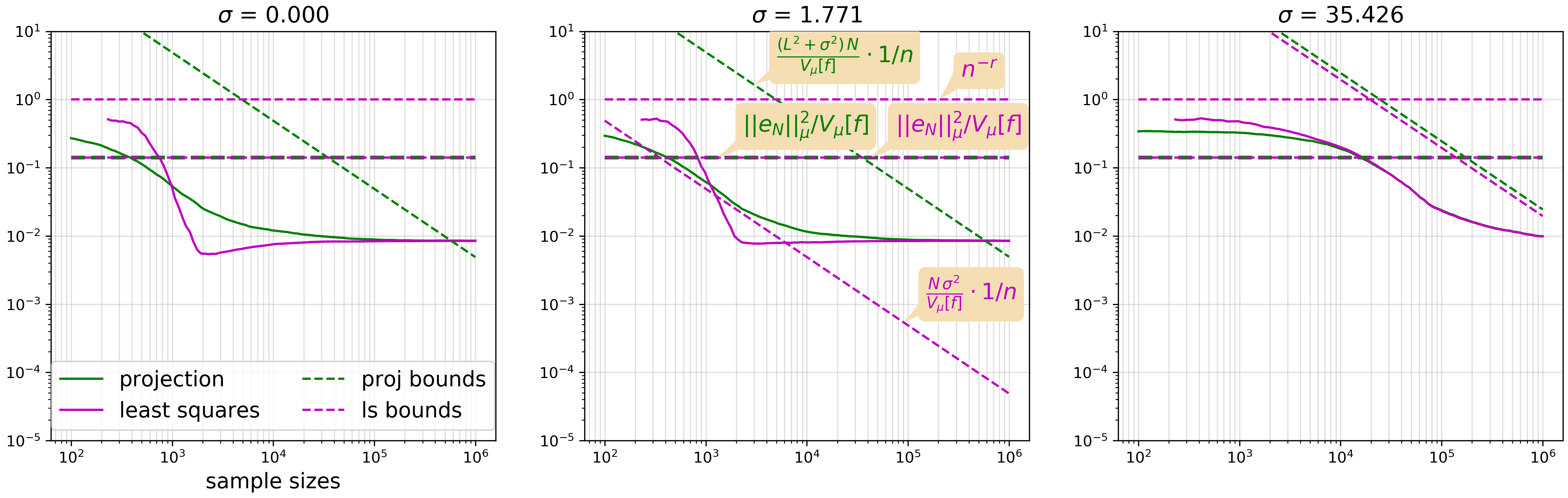}}
\caption{Ishigami function. Quadratic risk of Sobol' and total-effect indices $\max_{\ub} \big\{\E (S_\ub - \hat{S}_\ub)^2, \; \E (T_\ub - \hat{T}_\ub )^2 \big\}$ and the components of bounds (\ref{eq: si_mse_int}, \ref{eq: si_mse_ls_noise}), depending on the regressors number and noise level.}
\end{figure}

In this section, we will provide an illustration for Theorem~\ref{th: si_mse_int} and \ref{th: si_mse_ls_noise}. We numerically estimate the quadratic risk of Sobol' and total-effect indices $\max_{\ub} \big\{\E (S_\ub - \hat{S}_\ub)^2, \; \E (T_\ub - \hat{T}_\ub )^2 \big\}$ depending on the sample sizes and show how it relates to the components of risk bounds (\ref{eq: si_mse_int}, \ref{eq: si_mse_ls_noise}) 
for the two methods of approximation construction, that is (up to additional constant factors)
\[
\frac{\|e_N\|_{\mu}^2 }{\V_{\mu}[f]}, \;\;
  \frac{L^2 + \sigma^2}{\V_{\mu}[f]}  \cdot \frac{N}{n}, \;\;
\frac{\sigma^2}{\V_{\mu}[f]}  \cdot \frac{N}{n} \;\;
\text{and} \;
n^{-r}.
\]
These components allow us to understand the behavior of the risk better than the obtained risk bounds, which can significantly overestimate the risk.

Numerical risk estimation is based on $n_{run} = 100$ independent random designs (and noise realizations if $\sigma > 0$). In order to estimate the best theoretical error $\|e_N\|_{\mu}^2$, we follow (\ref{eq: MSE_noiseless_app}) and the least squares method using training and test samples of size $10^6$ each. Following (\ref{eq: def_kappa_r_and_condition_KN_n}), $r$ is calculated as
\begin{equation}
\label{eq: r_calculaton}
r = \max \, \left\{
0,  \;\;
\frac{3 \cdot \ln(3/2) - 1}{2} \cdot \frac{n}{K_N \ln n} - 1 
\right\}.
\end{equation}

{\bf Approximating function.} We use approximations of the test functions based on multivariate Legendre  (\ref{eq: legendre_multivariate}) and Trigonometric (\ref{eq: Trigonometric}) polynomials, the least squares and the projection method, and the maximum degree truncation scheme (\ref{eq: max_degree}).

{\bf The effects of regressors number and the type of regressors} are shown in Figures~\ref{fig: Sobol_risk_ls}, \ref{fig: Sobol_risk_proj}, \ref{fig: Ishigami_risk_ls} and \ref{fig: Ishigami_risk_proj}. Observations are noiseless. For comparison, we use regressors based on Legendre and Trigonometric polynomials, which are both orthogonal w.r.t. a uniform measure. The parameters of the truncation schemes are the following:
\begin{itemize}
\item {Sobol' g-function:} $\alpha_{\max} \in \{4, 6, 9\}$, $N \in \{25, 49, 100\}$ regressors;
\item {Ishigami function:}   $\alpha_{\max} \in \{4, 5, 6\}$, $N \in \{125, 216, 343\}$ regressors;
\end{itemize}

The experiments demonstrate that the best theoretical error $\|e_N\|_{\mu}^2$  does influence ``limiting risk'', which corresponds to almost constant values of estimated risk at large sample sizes. Although this best error can be much bigger than the indices risk, larger value of $\|e_N\|_{\mu}^2$ in most cases leads to larger limiting risk.  As expected, an increase in the number of regressors $N$ corresponds to a decrease in limiting risk but also less stable estimates for small sample sizes. Note that different $N$ may correspond to almost the same approximation errors for the specific functions: compare $N=125$ and $N=216$ regressors in Figure~\ref{fig: Ishigami_risk_ls}. In addition, the truncation scheme has a  strong influence on the convergence speed (see, for comparison, Figure~\ref{fig: theoretical_max_and_sum_bounds_Ishigami}).

In Figures~\ref{fig: Sobol_risk_ls} and \ref{fig: Ishigami_risk_ls}, one can see a big difference between the two regressors types when using least squares. Given the maximum degree truncation scheme, trigonometric polynomials with $K_N = N$ provide better stability for small sample sizes compared to Legendre polynomials with $K_N = N^2$. This stability is especially evident when limiting risks for the two types of regressors are close. The effect is not apparent in the case of the projection method. In general, regressors based on trigonometric polynomials perform better in all cases presented. Note, however, that Sobol' g-function and Ishigami can be extended correspondingly to $1$- and $2\pi$-periodic functions in each of their arguments, which provides better theoretical error $\|e_N\|_{\mu}^2$ for the  trigonometric approximation of these specific functions compared to ``aperiodic'' ones.

According to (\ref{eq: r_calculaton}), the minimal sample size that provides sensible estimates based on Theorem~\ref{th: si_mse_ls_noise} increases considerably with increasing dimension and the maximum degree of the truncation scheme. In particular, in the case of Ishigami function and $N = 125$ regressors of Legendre type, it needs at least $n = 2{,}102{,}432 > 10^6$ design points to achieve $r > 0$. 

Finally, we can briefly discuss how to select an approximating metamodel for the evaluation of Sobol' indices; particularly, how to choose its structure and the number of regressors?

From a theoretical point of view, the structure of the  metamodel should be adapted to the problem. In particular, according to (\ref{eq: N_opt_ls_leg_noise-free}, \ref{eq: si_ols_asymptotic_cheb}, \ref{eq: si_ols_asymptotic_trig}, \ref{eq: N_opt_proj_noise}, \ref{eq: N_opt_ls_noise}), the number of regressors should depend on the smoothness of the analyzed function and the size of the experimental design. Despite the fact that these results are theoretically justified, they are essentially asymptotic and require {\it a-priori} information on the analyzed function.  

From a practical perspective, it is important that,  according to Theorem~\ref{th: si_error_general}, more accurate (in the RMSE sense) metamodels are preferable for the estimation of Sobol' indices. Thus, one can customize the metamodel structure for global sensitivity analysis using standard methods of {\it model selection} and {\it feature selection} \cite{ElstatHastie2009}, which are based on cross-validation, penalization of the metamodel complexity \cite{Laurent2000}, an adaptive construction of the approximation \cite{BlatmanSudret10b,MR2060166}, etc.

{\bf The effect of noise in observations} is presented in Figures~\ref{fig: Sobol_risk_noise} and \ref{fig: Ishigami_risk_noise}. Independent Gaussian noise is added to the output of the test functions. The~standard deviation of the noise equals $0$, $L/10$ and $2 L$ with $L = \max_{\xb \in \XX} \big|f(\xb)\big|$. Regressors are based on Legendre polynomials. The parameters of the maximum degree truncation scheme are the following:
\begin{itemize}
\item {Sobol' g-function:} $\alpha_{\max}=4$, $N=25$ regressors;
\item {Ishigami function:}  $\alpha_{\max}=5$, $N=216$ regressors;
\end{itemize}

As one can see, in comparison with least squares, the projection method is more stable for very small sample sizes $n \sim N$. Given the moderate sample size and low  or zero noise level, the quadratic risk of the estimates for the least squares method decays significantly faster than $n^{-1}$ and correspondingly than the risk for the projection method. Large noise levels lead to similar behavior of the two methods, and the components of their risks bounds proportional to $N/n$ are close.

We can hypothesize that the three components of the obtained risk (\ref{eq: si_mse_ls_noise}), namely $n^{-r}$, $\sigma^2 / \V_{\mu}[f]  \cdot N / n$ and $\|e_N\|_{\mu}^2 / \V_{\mu}[f]$, correspond to three ranges of sample sizes with different rates of risk decay for the least squares method.

\section{Conclusion}

We address the problem of the accuracy and risk of metamodels-based Sobol' indices in the random design setting. We obtained a general relation between the accuracy of arbitrary metamodel and the error of estimated Sobol' indices; and the specific asymptotic and non-asymptotic relations for the two methods of parameters estimation for metamodels with tensor structure, including approximations based on multivariate Legendre, Chebyshev,  and Trigonometric polynomials. Based on the obtained relation, we propose a method for Sobol' indices quality control that demonstrates a good performance compared to an existing approach. 

The results indicate that the risk convergence with sample growth for the least squares method in the noiseless case can be much faster in comparison with the noisy case and with the projection method. In particular, it is demonstrated theoretically and (for a certain range of sample sizes $n$) experimentally that the decay of quadratic risk can be much faster than $n^{-1}$ with increasing sample size and the number of regressors (see Table~\ref{tab: SI_risk_all}). 

It is shown that the key factors that provide the possibility for fast convergence of Sobol' indices estimates are the absence of noise in the output of the analyzed function, its high smoothness, and low dimension. Besides, the estimates of small and close to $1$ indices tend to converge faster, and the obtained bounds and the proposed method are more accurate for such indices.



\bibliographystyle{imsart-number}

\bibliography{sample}

%

\newpage


\appendix

\section{Proofs}
\label{sec: Proofs}

\subsection{Proof of Theorem~\ref{th: si_error_general}}

1) We will prove the theorem not only for Sobol' indices and total-effect indices but also for a general type of indices that takes into account the influence of arbitrary subsets of variables groups. Starting from Sobol-Hoeffding decomposition, divide the partial variances into two arbitrary disjoint groups:
\begin{equation}
\label{eq: partial_var_hoeffding_subgroups}
\V_{\mu}[f(\xb)] 
= 
\sum_{\substack{\ub \subseteq \{1, \ldots, d\}}  } \V_{\mu}[f_\ub(\xb_\ub)] 
= 
\sum_{\ub \in \UU} \V_{\mu}[f_\ub(\xb_\ub)] 
+
\sum_{\ub \in \unUU} \V_{\mu}[f_\ub(\xb_\ub)] ,
\end{equation} 
where $\UU$ and $\unUU$ are disjoint subsets of the set of all subsets of $\{1, \ldots, d\}$, {\it i.e.} $\UU \cap \unUU = \emptyset$ and $\UU \cup \unUU = \big\{\ub\colon \ub \subseteq  \{1, \ldots, d\}\big\}$. Define {\it a general Sobol' index} of subset $\UU$ as
\begin{equation}
\label{eq: SI_general_subset}
G_{\UU} 
\eq
\frac{\sum_{\ub \in \UU} \V_{\mu}[f_\ub(\xb_\ub)] }{\V_{\mu}[f]}
=
\sum_{\ub \in \UU} S_{\ub}.
\end{equation}
 
Similarly, denote $\hat{G}_\UU \eq \sum_{\ub \in \UU} \hat{S}_{\ub}$. Note that $\UU = \big\{\ub \big\}$ corresponds to Sobol' index: $G_{\UU} = S_{\ub}$, and $\UU = \big\{\vb\colon \vb \subseteq  \{1, \ldots, d\}, \; \ub \cap \vb \neq \emptyset \big\}$ corresponds to total-effect index: $G_{\UU} = T_{\ub}$.

2) Notice that all Sobol' indices remain the same if we replace $f$ and $\hat{f}$ with the functions $g = a_1 \cdot f + a_2$ and $h = b_1 \cdot \hat{f} + b_2$ correspondingly, where $a_1,a_2, b_1,b_2 \in \R$,\;$a_1 \neq 0$, $b_1 \neq 0$. 

Consider the inequalities
\begin{equation}
\label{eq: min_norm_RHS}
\begin{aligned}
\frac{{\|f - \hat{f}\|^2_{\mu} }}{\V_{\mu}[f]}
&\geq
\frac{{\|(f-\E_{\mu}f) - (\hat{f}-\E_{\mu}\hat f)\|^2_{\mu} }}{\V_{\mu}[f]}
=
\left\| \frac{f-\E_{\mu}f}{\V_{\mu}^{1/2}[f]} - \frac{\hat{f}-\E_{\mu}\hat f}{\V_{\mu}^{1/2}[f]} \right\|^2_{\mu}
\\
&\geq
\min_{k \in \R}\left\| \frac{f-\E_{\mu}f}{\V_{\mu}^{1/2}[f]} - k \cdot \frac{\hat{f}-\E_{\mu}\hat f}{\V_{\mu}^{1/2}[f]} \right\|^2_{\mu}
\\
&=
\left\| \frac{f-\E_{\mu}f}{\V_{\mu}^{1/2}[f]} - \frac{\la f-\E_{\mu}f, \hat{f}-\E_{\mu}\hat f \ra_{\mu} }{\V_{\mu}[\hat f]} \cdot \frac{\hat{f}-\E_{\mu}\hat f}{\V_{\mu}^{1/2}[f]} \right\|^2_{\mu}.
\end{aligned}
\end{equation}

Note that $k_{min} = \argmin_{k\in \R}\|p - k \cdot q\|^2_{\mu}$ for $p, q \in L^2(\XX, \mu)$ corresponds to orthogonal projection of $p$ onto $q$ that leads to $\la q, p - k_{min} \cdot q\ra_{\mu} = 0$.

3) If $\la f-\E_{\mu}f, \, \hat{f}-\E_{\mu}\hat f \ra_{\mu} = 0$, then $\|f - \hat{f}\|^2_{\mu} \geq \|(f-\E_{\mu}f) - (\hat{f}-\E_{\mu}\hat f)\|^2_{\mu} = \V_{\mu}[f] + \V_{\mu}[\hat f] > \V_{\mu}[f]$ and (\ref{eq: si_error_general}) holds true, as its LHS does not exceed $1$. Besides, (\ref{eq: si_error_general_each_S}) is  also true, as
\begin{eqnarray*}
  \Big|G_\UU - \hat{G}_\UU \Big|
& = &
\left\{ \sqrt{G_\UU(1-\hat{G}_\UU)} + \sqrt{\hat{G}_\UU(1-G_\UU)}\right\} 
\\
&\times& 
\left| \sqrt{G_\UU(1-\hat{G}_\UU)} - \sqrt{\hat{G}_\UU(1-G_\UU)}\right| 
 \\
& \leq &
\sqrt{G_\UU(1-\hat{G}_\UU)} + \sqrt{\hat{G}_\UU(1-G_\UU)}.
\end{eqnarray*}

 Let $\cov_\mu[f, \hat f] \eq \la f-\E_{\mu}f, \, \hat{f}-\E_{\mu}\hat f \ra_{\mu}  \neq 0$. Following (\ref{eq: min_norm_RHS}), define $g, h \in L^2(\XX, \mu)$ as 
\begin{equation}
\label{eq: def_g_h}
g \eq \frac{f-\E_{\mu}f}{\V_{\mu}^{1/2}[f]}, 
\;\;
h \eq \frac{\cov_\mu[f, \hat f]}{\V^{1/2}_{\mu}[f] \cdot \V^{1/2}_{\mu}[\hat f]} \cdot \frac{\hat{f}-\E_{\mu}\hat f}{\V_{\mu}^{1/2}[\hat f]} .
\end{equation}

The functions $g$ and $h$ have the same Sobol' indices as $f$ and $\hat f$ correspondingly. In addition, 
\begin{equation}
\label{eq: lemma_orthog_condition}
\V_{\mu}[g] = 1, 
\;\;
\E_{\mu}[g] = 0, 
\;\;
0 < \V_{\mu}[h] \leq 1,
\;\;
\E_{\mu}[h] = 0, 
\;\;
\la h, g-h\ra_{\mu} = 0,
\end{equation}
and
\begin{equation}
\label{eq: norm_f_fatf_g_h}
{\|g - h\|^2_{\mu}} = \frac{{\|g - h\|^2_{\mu} }}{\V_{\mu}[g]}
\leq
\frac{{\|f - \hat{f}\|^2_{\mu} }}{\V_{\mu}[f]}.
\end{equation}

Thus, it is sufficient to prove the theorem for the functions $g$ and $h$. To prove that, we will show that for any $\UU$ defined above
\begin{equation}
\label{eq: si_error_general_each_S_g_h}
\Big|G_\UU - \hat{G}_\UU \Big|
\leq  
\left\{ \sqrt{G_\UU(1-\hat{G}_\UU)} + \sqrt{\hat{G}_\UU(1-G_\UU)}\right\} 
\cdot 
{\|g - h\|_{\mu}}.
\end{equation}

4) Using Sobol-Hoeffding decompositions of $g$ and $h$
\[
g(\xb) 
=  
\sum_{\ub \subseteq \{1, \ldots, d\}} g_\ub(\xb_\ub) = g_\UU + g_{\unUU},
\;\;
h(\xb) 
=  
\sum_{\ub \subseteq \{1, \ldots, d\}} h_\ub(\xb_\ub)= h_\UU + h_{\unUU},
\]
where $g_\UU \eq \sum_{\ub \in \UU} g_\ub(\xb_\ub) $, $g_{\unUU} \eq \sum_{\ub \in \unUU} g_\ub(\xb_\ub)$, $h_\UU \eq \sum_{\ub \in \UU} h_\ub(\xb_\ub)$, $h_{\unUU} \eq \sum_{\ub \in \unUU} h_\ub(\xb_\ub)$. In addition, based on the decomposition properties
\[
\la g_\UU, g_{\unUU} \ra_{\mu} = 0, 
\;\;
\la h_\UU, h_{\unUU} \ra_{\mu} = 0,
\;\;
\la g_\UU, h_{\unUU} \ra_{\mu} = 0, 
\;\;
\la h_\UU, g_{\unUU} \ra_{\mu} = 0.
\]
We have the following representations of general Sobol' indices (\ref{eq: SI_general_subset})
\begin{equation}
\label{eq: SI_general_g_h}
G_\UU 
=
\frac{\|g_{\UU}\|^2_{\mu}}{\|g\|^2_{\mu}},
\;\;
\hat{G}_\UU 
=
\frac{\|h_{\UU}\|^2_{\mu}}{\|h\|^2_{\mu}}.
\end{equation}

5) We will need an auxiliary bound based on (\ref{eq: lemma_orthog_condition}) and Cauchy-Schwarz inequality:
\[
\|h\|^2_{\mu} 
= 
\la g, h\ra_{\mu}
=
\la g_\UU, h_\UU\ra_{\mu} + \la g_{\unUU}, h_{\unUU}\ra_{\mu}
\leq
\|g_{\UU}\|_{\mu} \|h_{\UU}\|_{\mu} + \|g_{\unUU}\|_{\mu} \|h_{\unUU}\|_{\mu},
\]
and therefore,
\begin{equation}
\label{eq: auxiliary_inequality_product}
\|h\|^4_{\mu} 
-
\|g_{\UU}\|^2_{\mu} \|h_{\UU}\|^2_{\mu} 
- 
\|g_{\unUU}\|^2_{\mu} \|h_{\unUU}\|^2_{\mu}
\leq 
2 \cdot \|g_{\UU}\|_{\mu} \|h_{\UU}\|_{\mu} \|g_{\unUU}\|_{\mu} \|h_{\unUU}\|_{\mu}.
\end{equation}

%

6) Consider LHS of (\ref{eq: si_error_general_each_S_g_h}). Denote $\alpha \eq \arccos G_{\UU}^{1/2} \in [0, \pi/2]$,  $\beta \eq \arccos \hat{G}_{\UU}^{1/2} \in [0, \pi/2]$ (see Fig.~\ref{fig:2d_example} for illustration in 2D case), then
\begin{eqnarray}
\Big|G_\UU - \hat{G}_\UU \Big| 
& = & 
\big|\cos^2 \alpha - \cos^2 \beta \big|
=
1/2 \cdot \big|\cos 2\alpha - \cos 2\beta \big|
\nonumber
\\
& = & 
\sin(\alpha+\beta) \cdot \big|\sin(\alpha - \beta) \big| 
\nonumber 
\\
\label{eq: S_diff_decomp_sin}
& = & 
\left\{ \sqrt{G_\UU(1-\hat{G}_\UU)} + \sqrt{\hat{G}_\UU(1-G_\UU)}\right\} 
\cdot
\big|\sin(\alpha - \beta) \big|.
\end{eqnarray}

7) Compare $\omega \eq \big|\sin(\alpha - \beta) \big|$ and $\|g - h\|_{\mu}$. Based on (\ref{eq: SI_general_g_h}, \ref{eq: auxiliary_inequality_product}),
\begin{equation}
\begin{aligned}
\label{eq: sin_vs_g-h}
\omega^2
&= 
\left\{ \sqrt{G_\UU(1-\hat{G}_\UU)} - \sqrt{\hat{G}_\UU(1-G_\UU)}\right\}^2
\\
&=
\left\{
\frac{\|g_{\UU}\|_{\mu} \|h_{\unUU}\|_{\mu}}{\|g\|_{\mu} \|h\|_{\mu}} 
-
\frac{\|g_{\unUU}\|_{\mu} \|h_{\UU}\|_{\mu}}{\|g\|_{\mu} \|h\|_{\mu}}
\right\}^2
\\
&=
\frac{1}{\|h\|^2_{\mu}} 
\Big\{
\|g_{\UU}\|^2_{\mu} \cdot \|h_{\unUU}\|^2_{\mu}
+
\|g_{\unUU}\|^2_{\mu} \cdot \|h_{\UU}\|^2_{\mu}
-
2 \cdot \|g_{\UU}\|_{\mu}  \|g_{\unUU}\|_{\mu} \|h_{\UU}\|_{\mu} \|h_{\unUU}\|_{\mu}
\Big\}
\\
&\leq
\frac{1}{\|h\|^2_{\mu}} 
\Big\{
\|g_{\UU}\|^2_{\mu}  \|h_{\unUU}\|^2_{\mu}
+
\|g_{\unUU}\|^2_{\mu}  \|h_{\UU}\|^2_{\mu}
-
\|h\|^4_{\mu} 
\\
&+
\|g_{\UU}\|^2_{\mu} \|h_{\UU}\|^2_{\mu} 
+ 
\|g_{\unUU}\|^2_{\mu} \|h_{\unUU}\|^2_{\mu}
\Big\}
\\
&= 
\frac{1}{\|h\|^2_{\mu}} 
\Big\{
\big(\|g_{\UU}\|^2_{\mu} + \|g_{\unUU}\|^2_{\mu} \big) 
\cdot
\big(\|h_{\UU}\|^2_{\mu} + \|h_{\unUU}\|^2_{\mu}\big)
-
\|h\|^4_{\mu} 
\Big\}
\\
&=
1 - \|h\|^2_{\mu} 
=
\|g - h\|^2_{\mu},
\end{aligned}
\end{equation}
where the last equality follows from (\ref{eq: lemma_orthog_condition}):
\[
\|g - h\|^2_{\mu} 
= 
\|g\|^2_{\mu} + \|h\|^2_{\mu} - 2 \la g, h \ra_{\mu} 
=
1 - \|h\|^2_{\mu}.
\]

8) Combining the results (\ref{eq: norm_f_fatf_g_h}, \ref{eq: S_diff_decomp_sin}, \ref{eq: sin_vs_g-h}),  obtain (\ref{eq: si_error_general_each_S_g_h}) and (\ref{eq: si_error_general_each_S}):
\begin{equation}
\label{eq: general_si_bound}
\begin{aligned}
\Big|G_\UU - \hat{G}_\UU \Big| 
&\leq
\left\{ \sqrt{G_\UU(1-\hat{G}_\UU)} + \sqrt{\hat{G}_\UU(1-G_\UU)}\right\} 
\cdot
\|g - h\|_{\mu}
\\
&\leq
\left\{ \sqrt{G_\UU(1-\hat{G}_\UU)} + \sqrt{\hat{G}_\UU(1-G_\UU)}\right\} 
\cdot
\frac{{\|f - \hat{f}\|_{\mu} }}{\V^{1/2}_{\mu}[f]}.
\end{aligned}
\end{equation}

9) Finally, as $\sqrt{G_\UU(1-\hat{G}_\UU)} + \sqrt{\hat{G}_\UU(1-G_\UU)} = \sin(\alpha + \beta) \leq 1$,  (\ref{eq: si_error_general}) holds true.  
\QEDA

\subsection{Proof of Corollary~\ref{th: si_error_general_resolved}}

1) The case $\mer = 0$ is trivial. Let $\mer \neq 0$. The proof for Sobol' indices and total-effects is the same.   Denote $\delta \eq \big|S_\ub - \hat{S}_\ub \big| \in [0, 1]$ and consider the consequence of Theorem~\ref{th: si_error_general} for small indices:
\begin{equation}
\label{eq: si_error_case_small}
\begin{aligned}
\delta
&\leq  
\left\{ \sqrt{S_\ub(1-\hat{S}_\ub)} + \sqrt{\hat{S}_\ub(1-S_\ub)}\right\} 
\cdot 
\mer
\\
&\leq
\Big\{S^{1/2}_\ub + \hat{S}^{1/2}_\ub \Big\}
\cdot 
\mer
\leq
\Big\{\sqrt{S_\ub} + \sqrt{S_{\ub}+\delta} \Big\}
\cdot 
\mer.
\end{aligned}
\end{equation}

2) The solution of the previous inequality w.r.t. $\delta$ is  the union of  intervals given as

\[
\frac{\delta^2}{\mer^2} - 2 S_{\ub} - \delta \leq 0
\;\;\;
\cup
\;\;\;
    \begin{cases}
      \frac{\delta^2}{\mer^2} - 2 S_{\ub} - \delta \geq 0, \\ \left(\frac{\delta^2}{\mer^2} - 2 S_{\ub} - \delta \right)^2 \leq 4 S_{\ub} \,(S_{\ub}+\delta).
    \end{cases}
\]
Based on these inequalities, obtain the bounds for $\delta$:
\[
\delta \leq 1/2 \cdot \big(\mer^2 + \sqrt{\mer^4 + 8 S_{\ub} \mer^2} \big)
\;\;\;
\cup
\;\;\;
\delta \leq \mer^2 + 2 \sqrt{S_{\ub}} \cdot \mer.
\]
Taking into account $1/2 \cdot \big(\mer^2 + \sqrt{\mer^4 + 8 S_{\ub} \mer^2} \big) \leq \mer^2 + \sqrt{2} \, \sqrt{S_{\ub}} \cdot \mer$ and $\delta \leq \mer$, we have
\begin{equation}
\label{eq: si_error_case_small_sol}
\Big|S_\ub - \hat{S}_\ub \Big|
\leq  
\min\big(1, \; 2 \,\sqrt{S_{\ub}} + \mer \big) \cdot \mer.
\end{equation}


3) Consider the consequence of Theorem~\ref{th: si_error_general}  for  the indices close to $1$:
\begin{eqnarray*}
\delta 
& \leq &  
\left\{ \sqrt{S_\ub(1-\hat{S}_\ub)} + \sqrt{\hat{S}_\ub(1-S_\ub)}\right\} 
\cdot 
\mer
\leq
\Big\{\sqrt{1-S_\ub} + \sqrt{1-\hat{S}_\ub} \Big\}
\cdot 
\mer
\\
& \leq &
\Big\{\sqrt{1 - S_\ub} + \sqrt{1-S_{\ub}+\delta} \Big\}
\cdot 
\mer,
\end{eqnarray*}
where we used $\big| (1-S_\ub) - (1 - \hat{S}_\ub)\big|
= \delta$.   The obtained inequality is similar to (\ref{eq: si_error_case_small}) and the corresponding solution is 
\begin{equation}
\label{eq: si_error_case_large_sol}
\Big|S_\ub - \hat{S}_\ub \Big|
\leq  
\min\big(1, \; 2 \,\sqrt{1-S_{\ub}} + \mer \big) \cdot \mer.
\end{equation}

4) The final result is the combination of (\ref{eq: si_error_case_small_sol}) and (\ref{eq: si_error_case_large_sol}).
\QEDA

\subsection{Proof of Corollary~\ref{th: si_error_general_max}}

1) Denote $\theta_{\ub} \eq \sign(S_\ub - \hat{S}_\ub) \in \{-1,0,1\}$ for $\ub \subseteq \{1, \ldots, d\}$. Define $\UU \eq \{\ub\colon \theta_{\ub} \geq 0\}$ and $\unUU \eq \{\ub\colon \theta_{\ub}=-1\}$. 

Using the definition (\ref{eq: SI_general_subset}) of general Sobol' index and the inequality (\ref{eq: general_si_bound}) from the proof of Theorem~\ref{th: si_error_general}, obtain 
\begin{equation}
\label{eq: sum_abs_si_bound}
\begin{aligned}
\sum_{\substack{\ub \subseteq \{1, \ldots, d\}}}\big|S_\ub - \hat{S}_\ub \big|
&= 
\sum_{\ub} \theta_{\ub} \cdot \big(S_\ub - \hat{S}_\ub \big)
=
G_{\UU} - \hat{G}_{\UU} + \hat{G}_{\unUU} - G_{\unUU}
\\
&=
\big|G_{\UU} - \hat{G}_{\UU}\big| + \big|\hat{G}_{\unUU} - G_{\unUU}\big|
\leq
2 \cdot \mer.
\end{aligned}
\end{equation}

2) Let  $\vb \subseteq \{1, \ldots, d\}$. Based on (\ref{eq: sum_abs_si_bound}), we have
\[
\Big\{\sum_{\ub}\big|S_\ub - \hat{S}_\ub \big| \Big\}^2
=
\sum_{\ub}\big(S_\ub - \hat{S}_\ub \big)^2
+
\sum_{\substack{\ub, \vb\\ \ub \neq \vb}}\big|S_\ub - \hat{S}_\ub \big| \cdot \big|S_\vb - \hat{S}_\vb \big|
\leq
4 \cdot \mer^2.
\]

On the other hand,
\[
0 = \Big\{\sum_{\ub}(S_\ub - \hat{S}_\ub) \Big\}^2
=
\sum_{\ub}\big(S_\ub - \hat{S}_\ub \big)^2
+
\sum_{\substack{\ub, \vb\\ \ub \neq \vb}}\big(S_\ub - \hat{S}_\ub \big) \cdot \big(S_\vb - \hat{S}_\vb \big).
\]

Hence
\[
\sum_{\ub}\big(S_\ub - \hat{S}_\ub \big)^2
=
- \sum_{\substack{\ub, \vb\\ \ub \neq \vb}}\big(S_\ub - \hat{S}_\ub \big) \cdot \big(S_\vb - \hat{S}_\vb \big)
\leq
\sum_{\substack{\ub, \vb\\ \ub \neq \vb}}\big|S_\ub - \hat{S}_\ub \big| \cdot \big|S_\vb - \hat{S}_\vb \big|.
\]

Finally, (\ref{eq: si_error_general_sum_squared}) follows from
\[
2 \cdot \sum_{\ub}\big(S_\ub - \hat{S}_\ub \big)^2
\leq
\sum_{\ub}\big(S_\ub - \hat{S}_\ub \big)^2
+
\sum_{\substack{\ub, \vb\\ \ub \neq \vb}}\big|S_\ub - \hat{S}_\ub \big| \cdot \big|S_\vb - \hat{S}_\vb \big|
\leq
4 \cdot \mer^2.
\]
\QEDA

\subsection{Proof of Theorem~\ref{th: si_error_general_strict}}

1) For any fixed subset of variables $\ub \subseteq \{1, \ldots, d\}$  consider multi-index $\al = (\alpha_1, \ldots ,\alpha_d) \in \NN^d$ such that $\alpha_i = 1$ if $i \in \ub$, and $\alpha_i = 0$ otherwise. Besides, let $\be = (\beta_1, \ldots ,\beta_d) \in \NN^d$ be the multi-index such that $\beta_i = 1 - \alpha_i$ for $i = 1, \ldots, d$. Due to the assumption given in Remark~\ref{rem: finite_set_discrete_measure}, there exists $\Psi_{\al}, \,\Psi_{\be} \in L^2(\XX, \mu)$ which are the elements of the function set (\ref{eq: ONB}).

2) Select two functions in $L^2(\XX, \mu)$:
\begin{equation*}
\begin{aligned}
f(\xb) 
&= 
\sqrt{S_\ub} \cdot \Psi_{\al}(\xb_\ub) +  \sqrt{1-S_\ub} \cdot \Psi_{\be}(\xb_{\unsim\ub}),
\\
\hat{f}(\xb) 
&= 
\delta_1 \cdot \sqrt{\hat{S}_\ub} \cdot \Psi_{\al}(\xb_\ub) +  \delta_1 \cdot \sqrt{1-\hat{S}_\ub} \cdot \Psi_{\be}(\xb_{\unsim\ub}),
\end{aligned}
\end{equation*} 
where $\delta_1 \eq \sqrt{S_\ub \cdot \hat{S}_\ub} + \sqrt{(1-S_\ub)(1- \hat{S}_\ub)} > 0$.  It is easy to see that these functions have  Sobol' (and total) indices equal $S_{\ub}$ and $\hat{S}_{\ub}$ for $\xb_{\ub}$ variables correspondingly. One can check that
\[
\mer^2 
=
1 - \delta_1^2
=
\left| \sqrt{S_{\ub}(1-\hat{S}_{\ub})} - \sqrt{\hat{S}_{\ub}(1-S_{\ub})}\right|^2.
\]

3) Hence $\Big\{ \sqrt{S_{\ub}(1-\hat{S}_{\ub})} + \sqrt{\hat{S}_{\ub}(1-S_{\ub})}\Big\} 
\cdot
\mer
=
\big|S_{\ub} - \hat{S}_{\ub} \big| 
$. 
\QEDA

\subsection{Proof of Theorem~\ref{th: si_error_general_lower_max}}
%

1) Following Section~\ref{sec: Tightness_of_bounds}, consider two functions in $L^2(\XX, \mu)$:
\begin{eqnarray*}
f(\xb) 
&=& 
\frac{1}{2} \cdot \Psi_{(1,0, \dots)}(\xb) +  \frac{\sqrt{3}}{2} \cdot \Psi_{(0,1, \dots)}(\xb),
\\
\hat{f}(\xb) 
&=& 
\frac{3}{4} \cdot \Psi_{(1,0, \dots)}(\xb) +  \frac{\sqrt{3}}{4} \cdot \Psi_{(0,1, \dots)}(\xb) + \delta_2 \cdot \Psi_{(0,0, \dots)}(\xb),
\end{eqnarray*}
where $\delta_2 \in \R$. According to (\ref{eq: example_strict_bound_degree30}, \ref{eq: example_bound_max_3}), if $\delta_2 = 0$, then (\ref{eq: si_error_general_lower_max}) holds true with $\eps=1$. To obtain (\ref{eq: si_error_general_lower_max}) for $\eps \in (0, 1]$, one can set $\delta_2 = 1/2 \cdot \sqrt{1-\eps^2}/\eps$. Indeed,  
\[
\eps \cdot \frac{{\|f - \hat{f}\|_{\mu} }}{\V^{1/2}_{\mu}[f]}
=
\eps \cdot [1/4 + \delta_2^2]^{1/2} 
= 
\frac{1}{2} 
=
\max_{\ub}\Big|S_\ub - \hat{S}_\ub \Big| 
=
\max_{\ub}\Big|T_\ub - \hat{T}_\ub \Big| 
.
\]

2) In the case of $\eps=0$ it is sufficient to consider the functions: $f(\xb) =  \Psi_{(1,0, \dots)}(\xb) $ and $\hat{f}(\xb) = 2 \cdot \Psi_{(1,0, \dots)}(\xb)$.
\QEDA

\subsection{Proof of Theorem~\ref{th: si_mae_risk_general}}


1) We use the notation $S_\ub$ for both Sobol' indices and total-effects. Following Condition~\ref{cond: random_design} and (\ref{eq:training_sample}),  define the measure $\rho$ over $\XX \times \R$ and the measure $d\brho(\des, Y) = \otimes_{i=1}^n d\rho(\xb_i, y_i)$ over the training samples space $\XX^n \times \R^n$.

%


2) Taking into account  (\ref{eq: appr_deterministic_training}), consider separately two domains in the training samples space $\XX^n \times \R^n$: $\chi^{0} \eq \{(\des, Y)\colon \; \V_{\mu} [\hat{f}]= 0 \}$ and $\chi^{+} \eq  \XX^n \times \R^n \setminus{\chi^{0}}$.

3) According to Theorem~\ref{th: si_error_general}, if $\V_{\mu}[\hat{f}] > 0$, then 
\begin{equation}
\label{eq: bound_S_squared_all_cases}
(S_\ub - \hat{S}_\ub \big)^2
\leq
\frac{{\|f - \hat{f}\|^2_{\mu} }}{\V_{\mu}[f]}.
\end{equation}
If $\V_{\mu}[\hat{f}] = 0$, then $\|f - \hat{f}\|_{\mu}^2 \geq \V_{\mu}[f]$, and (\ref{eq: bound_S_squared_all_cases}) also holds true, taking into account Sobol' indices redefined for constant approximations.


4) The risk takes the form:
\begin{equation}
\label{eq: si_mse_risk_general_parts}
\begin{aligned}
\E\big(S_\ub - \hat{S}_\ub \big)^2  
&= 
\left\{
\int_{\chi^{+}}+\int_{\chi^{0}} 
\right\}
\big(S_\ub - \hat{S}_\ub \big)^2 d\brho(\des, Y) 
\\
&\leq
\frac{{\E \|f - \hat{f}\|_{\mu}^2 }}{\V_{\mu}[f]}
\eq
\ris^2,
\end{aligned}
\end{equation}
where we used (\ref{eq: bound_S_squared_all_cases}).

5) Note that $\mer^2 + 2\sqrt{S_{\ub}} \, \mer  \geq \mer^2 \geq 1$ if $\V_{\mu}[\hat{f}] = 0$. Similarly to (\ref{eq: si_mse_risk_general_parts}), we have based on Corollary~\ref{th: si_error_general_resolved}
\begin{eqnarray*}
\E\big|S_\ub - \hat{S}_\ub \big| 
&\leq& 
\E \left[\mer^2 + 2\sqrt{S_{\ub}} \, \mer \right] 
\\
&\leq& 
\E \mer^2 + 2 \sqrt{S_{\ub}} \, \sqrt{\E \mer^2}
= \ris^2 + 2 \sqrt{S_{\ub}} \cdot \ris.
\end{eqnarray*}
\QEDA

\subsection{Proof of Corollary~\ref{th: si_sum_mse_risk_general}}
The proof repeats the proof of Theorem~\ref{th: si_mae_risk_general} using Corollary~\ref{th: si_error_general_max} and taking into account that if $\V_{\mu} [\hat{f}] = 0$, then for $\ub \subseteq \{1, \ldots, d\}$
\begin{eqnarray*}
\sum_{\ub} \big(S_\ub - \hat{S}_\ub \big)^2 
&=&
\sum_{\ub} \big(S_\ub - 2^{-d} \big)^2 
=
\sum_{\ub} S_\ub^2 - 2 \cdot 2^{-d} \cdot 1 + 2^{-2d} \cdot 2^{d}
\\
&=&
\sum_{\ub} S_\ub^2 - 2^{-d}
<
\Big\{\sum_{\ub} S_\ub \Big\}^2
=
1
<
2
\leq
2 \cdot \mer^2.
\end{eqnarray*}
\QEDA

\subsection{Proof of Theorem~\ref{th: si_mse_int}}

1) By definition (\ref{eq: best_apprx_pnu_N}), the theoretical  orthogonal projection $f_N$ is a linear combination of $\{\Psi_{\al}, \; \al \in \rr_N\}$:
\[
f_N(\xb) =  \sum_{\al \in \rr_N} c_{\al} \Psi_{\al}(\xb),
\]
with some $c_{\al} \in \R$.

2) Let $\hat{f}^{P}$ be  an approximation based on the training sample of size $n$ and $\{\hat{c}_{\al}\}$ be the corresponding estimated coefficients. For $\hat{c}_{\al}$ defined as (\ref{eq: Projection_coeff}) it holds for the expectation
\[
\E [\hat{c}_{\al}] = \E \big[(f(\xb)+\eta) \Psi_{\al}(\xb) \big] 
=
\E \big[ f(\xb) \Psi_{\al}(\xb)\big] 
=
c_{\al},
\]
and for the variance
\begin{eqnarray*}
\V [\hat{c}_{\al}] 
&=& 
\frac{1}{n} \V  \big[(f(\xb)+\eta)\Psi_{\al}(\xb)\big]
\leq
\frac{1}{n} \E  \big[(f(\xb)+\eta)^2 \Psi^2_{\al}(\xb)\big]
\\
&=&
\frac{1}{n} \E  \big[f^2(\xb) \Psi^2_{\al}(\xb)\big] + \frac{1}{n}  \E  \big[\eta^2 \Psi^2_{\al}(\xb)\big]
\leq
\frac{1}{n} (L^2 + \sigma^2).
\end{eqnarray*}

3)
We have $\|f - \hat{f}^{P}\|_{\mu}^2 
= 
\|f - f_N\|_{\mu}^2 +  \|f_N - \hat{f}^{P}\|_{\mu}^2 
=
\|e_N\|_{\mu}^2 +  \|f_N - \hat{f}^{P}\|_{\mu}^2 
=
\|e_N\|_{\mu}^2 +  \sum_{\al \in \rr_N} (c_{\al} - \hat{c}_{\al})^2$. Then,
\begin{eqnarray*}
\E \|f - \hat{f}^{P}\|_{\mu}^2 
&=&
\|e_N\|_{\mu}^2 +  \E   \sum_{\al \in \rr_N} (c_{\al} - \hat{c}_{\al})^2  
=
\|e_N\|_{\mu}^2 + \sum_{\al \in  \rr_N} \V[ \hat{c}_{\al}]
\\
&\leq&
\|e_N\|_{\mu}^2 + \frac{N}{n}(L^2+\sigma^2).
\end{eqnarray*}

4) It remains to apply Theorem~\ref{th: si_mae_risk_general}.
\QEDA

\subsection{Proof of Theorem~\ref{th: si_mse_ls_noise}}

\begin{lemma}
\label{le: omp_mse_ls_risk}
Let $\bchi$ be the subset of $\XX^n \times \R^n$ that includes training samples for which the spectral norm $|||\Phi^T \Phi/n - I_N||| \leq 1/2$. Under conditions of Theorem~\ref{th: si_mse_ls_noise} it holds
\[
\int_{\bchi} \|f - \hat{f}^{LS}\|_{\mu}^2  d\brho(\des, Y)
\leq
\left(1 + \frac{4\kappa_r}{\ln n}\right)\|e_N\|^2_{\mu} + 4\sigma^2 \cdot \frac{N}{n}.
\]
\end{lemma}

{\bf Proof of Lemma~\ref{le: omp_mse_ls_risk}}

1) The proof includes modified versions of Theorem~$2$ and $3$ in \cite{Cohen13}. Denote the expansion coefficients of the best approximation (\ref{eq: best_apprx_pnu_N}) as $\cb \in \R^N$. We have $f_N(\xb) = \cb^T \bPsi(\xb)$. Provided $\det (\Phi^T \Phi)~\neq~0$,
\begin{eqnarray*}
\hat{\cb} 
&=& 
(\Phi^T \Phi)^{-1}\Phi^T Y 
=
(\Phi^T \Phi)^{-1}\Phi^T \big(\Phi \cb + e_N(\des) + \etab \big)
\\
&=&
\cb + (\Phi^T \Phi)^{-1}\Phi^T \big(e_N(\des)+ \etab \big),
\end{eqnarray*}
where $e_N(\des) \eq (e_N(\xb_1), \ldots, e_N(\xb_n))^T \in \R^n$. 

Hence
\[
 \|\cb - \hat{\cb}\| 
 =
\left\| \left(\frac{\Phi^T \Phi}{n} \right)^{-1}\frac{\Phi^T \big(e_N(\des)+ \etab \big)}{n} \right\|.
\]

2) By definition, $\bchi$  includes samples for which the spectral norm 
\begin{equation}
\label{eq: cond_design_ls}
|||\Phi^T \Phi/n - I_N||| \leq 1/2.
\end{equation}

From (\ref{eq: cond_design_ls}) one can obtain on $\bchi$ domain
\[
1/2 \leq |||\Phi^T \Phi/n||| \leq 3/2,
\]
\[
2/3 \leq |||(\Phi^T \Phi/n)^{-1}||| \leq 2.
\]

Then, provided $(\des, Y) \in \bchi$
\[
 \|\cb - \hat{\cb}\| 
\leq
2 \left\| \frac{\Phi^T \big(e_N(\des)+ \etab \big)}{n} \right\|.
\]

3) Back to the lemma's statement
\begin{eqnarray*}
\int_{\bchi} \|f - \hat{f}^{LS}\|_{\mu}^2  d\brho(\des, Y)
&\leq&
 \|e_N\|^2_{\mu} + \int_{\bchi}  \|\cb - \hat{\cb}\|^2  d\brho(\des, Y)
\\
&\leq&
\|e_N\|^2_{\mu} + \int_{\bchi} 4\left\| \frac{\Phi^T  \big(e_N(\des)+ \etab \big)}{n} \right\|^2 d\brho(\des, Y)
\\
&\leq&
\|e_N\|^2_{\mu} + 4 \, \E \left\| \frac{\Phi^T  \big(e_N(\des)+ \etab \big)}{n} \right\|^2
\\
&=&
\|e_N\|^2_{\mu} 
+ 
4 \, \E \left\| \frac{\Phi^T  e_N(\des)}{n} \right\|^2
+ 
4 \, \E \left\| \frac{\Phi^T \etab }{n} \right\|^2.
\end{eqnarray*}
The last equality is valid due to $\E \etab = 0$ and the independence of $\des$ and $\etab$.

4) Taking into account $\E [\Psi_{\al}(\xb_i) e_N(\xb_i)] = \E [ \Psi_{\al}(\xb_j) e_N(\xb_j)] = 0 $, we have
\begin{eqnarray*}
\E \left\| \frac{\Phi^T e_N(\des)}{n} \right\|^2 
&=&
\E \sum_{\al \in \rr_N} \left(\frac{1}{n}\sum_{i=1}^n\Psi_{\al}(\xb_i) e_N(\xb_i) \right)^2
\\
&=&
\frac{1}{n^2}  \E \sum_{\al \in \rr_N} \sum_{i=1}^n \sum_{j=1}^n \Psi_{\al}(\xb_i) e_N(\xb_i) \Psi_{\al}(\xb_j) e_N(\xb_j)
\\
&=&
\frac{1}{n^2}  \E \sum_{\al \in \rr_N} \sum_{i=1}^n \Psi^2_{\al}(\xb_i) e_N^2(\xb_i)  
=
\frac{1}{n}  \E_{\mu} \sum_{\al \in \rr_N} \Psi^2_{\al}(\xb) e_N^2(\xb)
\\
&\leq&
\frac{K_N}{n}  \E_{\mu}  [e_N^2(\xb)]
\leq
\frac{\kappa_r}{\ln n}  \|e_N\|^2_{\mu}.
\end{eqnarray*}

5) Similarly, using $\E [\Psi_{\al}(\xb_i) \cdot \eta_j] = 0 $,
\[
\E \left\| \frac{\Phi^T  \etab}{n} \right\|^2 
=
\E \sum_{\al \in \rr_N} \left(\frac{1}{n}\sum_{i=1}^n\Psi_{\al}(\xb_i) \cdot  \eta_i \right)^2
=
\frac{1}{n}  \E \sum_{\al \in \rr_N} \Psi^2_{\al}(\xb) \cdot \eta^2
=
\frac{N}{n} \sigma^2.
\]

6) Thus,
\[
\int_{\bchi} \|f - \hat{f}^{LS}\|_{\mu}^2  d\brho(\des, Y)
\leq
\left( 1 + \frac{4\kappa_r}{\ln n} \right) \|e_N\|^2_{\mu} + 4\sigma^2 \cdot \frac{N}{n} .
\]
\QEDA

{\bf Proof of Theorem~\ref{th: si_mse_ls_noise}}

1) We use the notation $S_\ub$ for both Sobol' indices and total-effects. Define two subsets of the training samples space  $\XX^n \times \R^n$. Let $\bchi \eq \{(\des, Y)\colon \; |||\Phi^T \Phi/n - I_N||| \leq 1/2 \}$ and $\chi^{0}$ be the subset of $\bchi$ that leads to constant approximation: $\chi^{0} \eq \{(\des, Y) \in \bchi\colon \; \V_{\mu} [\hat{f}^{LS}]= 0 \}$. Note that $|||\Phi^T \Phi/n - I_N|||$ depends only on the design $\des$ and is independent from $Y$.

2) The risk takes the form:
\begin{eqnarray*}
\E\big(S_\ub - \hat{S}^{LS}_\ub \big)^2
&=&
\left\{
\int_{\XX^n \times \R^n \backslash \bchi}  
+
\int_{\bchi \backslash \chi^{0}} 
+ \int_{\chi^{0}}
\right\}
\big(S_\ub - \hat{S}^{LS}_\ub \big)^2 d\brho(\des, Y)
\\
&\leq&
P\left\{\XX^n \times \R^n \backslash \bchi \right\} 
+ 
\left\{
\int_{\bchi \backslash \chi^{0}} 
+ 
\int_{\chi^{0}}
 \right\}
\big(S_\ub - \hat{S}^{LS}_\ub \big)^2 d\brho(\des, Y),
\end{eqnarray*} 
where we used $(S_\ub - \hat{S}^{LS}_\ub)^2 \leq 1$.

3) According to Theorem~\ref{th: si_error_general}, if $\V_{\mu}[\hat{f}^{LS}] > 0$, then 
\begin{equation}
\label{eq: bound_S_squared_ls_all}
(S_\ub - \hat{S}^{LS}_\ub \big)^2 \leq \mer^2.
\end{equation}
Since $\|f - \hat{f}^{LS}\|_{\mu}^2 \geq \V_{\mu}[f]$ provided $\V_{\mu}[\hat{f}^{LS}] = 0$, (\ref{eq: bound_S_squared_ls_all}) also holds true for the constant approximation. 



4) Using (\ref{eq: condition_K_N_and_n_05}, \ref{eq: bound_S_squared_ls_all}), obtain
\begin{eqnarray*}
\E\big(S_\ub - \hat{S}^{LS}_\ub \big)^2
&\leq&
2n^{-r}
+
\left\{
\int_{\bchi \backslash \chi^{0}} 
+ 
\int_{\chi^{0}}
 \right\}
\mer^2 d\brho(\des, Y)
\\
&\leq&
\int_{\bchi} \mer^2  d\brho(\des, Y)  + 2n^{-r}.
\end{eqnarray*}

5) Applying Lemma~\ref{le: omp_mse_ls_risk},
\[
\E\big(S_\ub - \hat{S}^{LS}_\ub \big)^2
\leq
\frac{\tau_n}{\V_{\mu}[f]}\|e_N\|^2_{\mu}
+ 
\frac{4 \sigma^2 }{\V_{\mu}[f]} \cdot \frac{N}{n} 
+
2n^{-r}
\leq
\ris_{LS}^2 + 2n^{-r},
\]
where $\tau_n \eq 1 + 4\kappa_r/\ln{n}$ with $\tau_n \to 1$ as $n \to \infty$. Besides, Condition~\ref{cond: ls_K_N_bound} of the theorem implies $\kappa_r \leq 1/2 \cdot [3 \cdot \ln(3/2) - 1]$ and $n \geq 33$. Then, corresponding values  $\tau_n < 1.2$.




6)  Similarly, based on Corollary~\ref{th: si_error_general_resolved}:
\begin{eqnarray*}
\E\big|S_\ub - \hat{S}^{LS}_\ub \big|
&\leq&
P\left\{\XX^n \times \R^n \backslash \bchi \right\} 
+ 
\left\{
\int_{\bchi \backslash \chi^{0}} 
+ 
\int_{\chi^{0}}
 \right\}
\big|S_\ub - \hat{S}^{LS}_\ub \big| d\brho(\des, Y)
\\
&\leq&
2n^{-r}
+
\left\{
\int_{\bchi \backslash \chi^{0}} 
+ 
\int_{\chi^{0}}
 \right\} \left(\mer^2 + 2\sqrt{S_{\ub}} \, \mer \right) d\brho(\des, Y)
\\
&=&
\int_{\bchi}
\left(\mer^2 + 2\sqrt{S_{\ub}} \, \mer \right) d\brho(\des, Y)
+
2n^{-r}
\\
&\leq&
\int_{\bchi} \mer^2 d\brho(\des, Y)
+
2\sqrt{S_{\ub}}  \cdot \left\{\int_{\bchi} \mer^2  d\brho(\des, Y)\right\}^{1/2}
+
2n^{-r}
\\
&\leq&
\ris_{LS}  \left(\ris_{LS} + 2\sqrt{S_{\ub}} \right) + 2n^{-r}. 
\end{eqnarray*}
\QEDA

\subsection{Proof of Corollary~\ref{th: converg_ls_asymp}}

The result follows from the inequality
\begin{eqnarray*}
\E\big(S_\ub - \hat{S}^{LS}_\ub \big)^2
&\leq&
\frac{1}{\V_{\mu}[f]} \left[ \left(1 + 4\frac{ K_N}{n}\right)\|e_N\|^2_{\mu} + 4\sigma^2 \cdot \frac{N}{n} \right]
\\
&+&
2 \exp \left[\ln N \cdot \left(1  - \frac{n}{K_N \ln N} \cdot c_{1/2} \right)\right],
\end{eqnarray*}
which is based on Lemma~\ref{le: spectral_norm_exponential_small} and the proofs of Theorem~\ref{th: si_mse_ls_noise} and Lemma~\ref{le: omp_mse_ls_risk}. 
\QEDA


\end{document}